\documentclass[USenglish]{article}	

\usepackage[utf8]{inputenc}				
\usepackage[big,online]{dgruyter}	
\usepackage{lmodern} 
\usepackage{microtype}
\usepackage[numbers,square,sort&compress]{natbib}

\usepackage{bm}
\usepackage{enumerate}
\graphicspath{{figs}}

\theoremstyle{dgthm}

\theoremstyle{dgdef}

\begin{document}

\title{Computational identification of the source domain in an inverse problem of potential theory}

\runningtitle{Identification of the source domain in an inverse potential problem}

\author*[1]{P.N.~Vabishchevich}
\runningauthor{P.N.~Vabishchevich}

\affil[1]{\protect\raggedright 
Lomonosov Moscow State University, 1, Leninskie Gory,  Moscow, Russia; 
North-Eastern Federal University, 58, Belinskogo st, Yakutsk, 677000, Russia, e-mail: vab@cs.msu.ru}
	
\dedication{the remarkable scientist and wonderful person, Professor Anatoly Yagola.}
	
\abstract{
The inverse potential problem consists in determining the density of the volume potential from measurements outside the sources.
Its ill-posedness is due both to the non-uniqueness of the solution and to the instability of the solution with respect to measurement errors.
The inverse problem is solved under additional assumptions about the sources using regularizing algorithms.
In this work, an inverse problem is posed for identifying the domain that contains the sources. 
The computational algorithm is based on approximating the volume potential by the single-layer potential on the boundary of the domain containing the sources.
The inverse problem is considered in the class of a priori constraints of nonnegativity of the potential density. 
Residual minimization in the class of nonnegative solutions is performed using the classical Nonnegative Least Squares algorithm.
The capabilities of the proposed approach are illustrated by numerical experiments for a two-dimensional test problem with an analytically prescribed potential on the observation surface.
}


\maketitle

\section{Introduction}

Inverse problems of potential theory \cite{isakov2006inverse} are of great practical importance, for example, in geophysics \cite{zhdanov2015inverse}.
In this case, it is necessary to reconstruct the characteristics of field sources from the known values of the field itself. 
The exterior inverse problem consists in determining the density of the volume (Newtonian) potential $\rho(\bm x)$ in the domain $D$ from the values of the potential $u(x)$ given on some surface $\Gamma$ lying outside the domain $D$.
Such problems are essentially ill-posed \cite{Tikhonov1977,Lavrentiev1986} due to the non-uniqueness of the solution and the instability of the solution to small changes in the input data.

The uniqueness of the solution to the inverse potential problem has been established under various restrictions on the desired potential density. For a given constant density, Novikov P.S. showed \cite{Novikov1938} that uniqueness holds in the class of star-shaped domains $D$. 
In the case of a given domain $D$, uniqueness is established, for example, for harmonic potential density, as well as if the density does not depend on one variable. These results have been generalized in various directions (see, e.g., \cite{isakov2006inverse,Isakov1990,Prilepko2000}).

In approximate solutions of inverse potential problems, various numerical methods are used.
For the nonlinear problem of finding the domain $D$ for a known density, the most promising approach is the use of the level-set method \cite{burger2001level}.
This approach is applied (see, e.g., \cite{hettlich1997recovery,lu2015local}) both directly to the volume potential and to the corresponding inverse problem for the Poisson equation.
For the linear inverse potential problem of identifying the right-hand side, special computational algorithms are used for the corresponding boundary value problems that take into account the features of the desired density \cite{samarskiui2007numerical}.

Since it is impossible to uniquely solve the inverse potential problem in the general case, one can pose more particular problems.
From a practical point of view, it would be interesting to localize the sources.
We distinguish some domain of potential source location $\Omega$. Based on the solution of an auxiliary inverse problem, we want to determine whether this domain contains all the sources (the domain $D$). Under these conditions, parametric calculations for different $\Omega$ would allow us to localize $D$.
In such a formulation, we are effectively continuing the solution toward the sources, from the observation surface $\Gamma$ to the boundary of the domain $\Omega$.
We attempt to answer the question of whether $D$ is localized within $\Omega$ by solving this continuation problem.

Various approaches can be used to compute the potential outside $\Omega$.
In this regard, we note the method of integral equations, which was applied, for example, in \cite{tikhonov1968prodolgenie} for solving the continuation problem. We work in the class of nonnegative densities, where $\rho(\bm x) \ge 0$ in the domain $D$. 
Under these conditions, the potential density of the single layer on $S$ ($D\subset \Omega$) is nonnegative \cite{vabishchevich2024computational}.
Such a priori constraints can be used in approximate solutions of the considered ill-posed potential problems \cite{vasin2013ill}. 

This work considers a particular inverse problem of assessing whether all sources $D$ belong to a certain given domain by approximating the potential values on $\Gamma$ with the potential values of a single layer supported on the boundary of the domain $\Omega$. The localization domain is identified based on the residual estimate for different $\Omega$.
The continuation problem for the potential is solved in the class of nonnegative densities.
The possibility of estimating the localization domain by the total mass of the single-layer potential is noted.
The capabilities of the proposed Nonnegative Density Domain algorithm are illustrated by numerical examples for a two-dimensional test problem.

\section{Problem statement}

\begin{figure}[ht]
\centering
\includegraphics[width=0.5\textwidth]{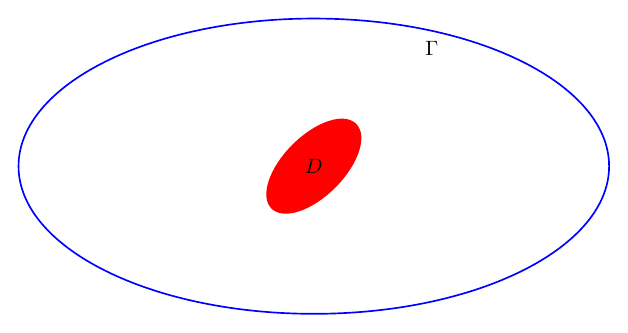}
\caption{To the problem statement.}
\label{f-1}
\end{figure}

We consider a two-dimensional problem of localizing the sources of the volume potential for the Laplace operator from data on the boundary of a domain that contains all sources.

The volume potential at the point $\bm{x}=(x^{(1)},x^{(2)})$ is defined by
\begin{equation}\label{2.1}
  u(\bm{x})=\int_{D}\rho(\bm{y})G(\bm{x},\bm{y})\, d \bm y,\qquad \bm{x}\in\mathbb{R}^{2} .
\end{equation}
Here the domain $D$ contains all sources, so that $\rho(\bm{y})=0$ for $\bm{y} \not\in D$.
For the kernel \(G(\bm{x},\bm{y})\) we have
\[
  G(\bm{x},\bm{y})= - \frac{1}{2 \pi} \ln{|\bm{x}-\bm{y}|},
 \quad |\bm{x}-\bm{y}|^{2}=\sum_{\alpha=1}^{2}(x^{(\alpha)}-y^{(\alpha)})^{2} .
\]

The inverse problem is posed as follows (Fig.~\ref{f-1}).
The potential is measured on the boundary $\Gamma$:
\begin{equation}\label{2.2}
  u(\bm{x}) = \varphi(\bm{x}),
  \quad \bm x \in \Gamma.
\end{equation}
In the general inverse potential problem, it is necessary to obtain information about the distribution of sources from these data.

In the above formulation, the solution of problem \eqref{2.1}, \eqref{2.2} for the identification of $\rho(\bm{y})$ is not unique. Therefore, this inverse potential problem is considered under additional assumptions.
In this regard, one can note the case when the density $\rho(\bm{y})$ is known and it is required to determine the domain $D$. 
A second class of inverse problems is associated with determining the spatial location of sources in a given domain $D$ under additional assumptions about $\rho(\bm{y})$ (for example, when the density depends only on one variable).

\begin{figure}[ht]
\centering
\includegraphics[width=0.5\textwidth]{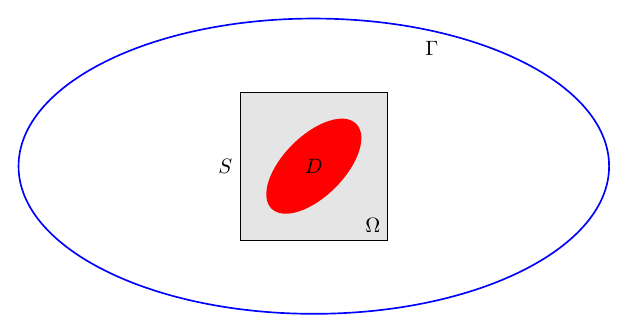}
\caption{Localization domain of the sources.}
\label{f-2}
\end{figure}

Within the framework of a more particular problem, we want to be able to identify whether the domain $D$ (Fig.~\ref{f-2}) lies entirely within a given domain $\Omega$ with boundary $S$, or not.
By performing computations with different positions of the domain $\Omega$, we would have the opportunity to localize the sources.

\section{Computational algorithms}

Computational identification of the source domain $\Omega$ is based on approximating the volume integral by a surface integral on the boundary $S$.
Various computational algorithms are applied for the numerical solution of the corresponding integral equation. 
Our approach is based on seeking an approximate solution in the class of nonnegative kernels of the Fredholm integral equation.

\subsection{Integral equation}   

The volume integral \eqref{2.1} outside the domain $\Omega$ is replaced by the surface single-layer potential
\begin{equation}\label{3.1}
  u(\bm x) = \int_{S}\mu(\bm{y})G(\bm{x},\bm{y})  \, d \bm y,
  \quad \bm x \notin \Omega \cup S .
\end{equation}

Determination of the single-layer potential density from the observation conditions \eqref{2.2} leads to a Fredholm integral equation of the first kind:
\begin{equation}\label{3.2}
\int_{S}\mu(\bm{y})G(\bm{x},\bm{y})  \, d \bm y = \varphi(\bm{x}), \quad \bm{x}\in\Gamma .
\end{equation}

For an approximate solution of the integral equation \eqref{3.2}, we can use various computational algorithms.
Taking into account the ill-posedness of this problem, we focus on regularization methods.

Let us write \eqref{3.2} in the form of an operator equation of the first kind
\begin{equation}\label{3.3}
  \mathcal{A} \mu = \varphi ,
\end{equation}
where
\[
  \mathcal{A} \mu = \int_{S}\mu(\bm{y})G(\bm{x},\bm{y})  \, d \bm y .
\]

To solve problem \eqref{3.3}, we can use the classical Tikhonov regularization method.
For a given regularization parameter $\alpha$, the approximate solution is defined from 
\begin{equation}\label{3.4}
  \|\mathcal{A} \mu - \varphi \|_\Gamma^2 + \alpha \|\mu\|_S^2 \rightarrow \min_{\mu} .
\end{equation}
In \eqref{3.4}, $\| \cdot\|_\Gamma, \ \| \cdot\|_S$ are the norms in $L_2(\Gamma)$ and $L_2(S)$:
\[
  \|\varphi\|_\Gamma^2 = \int_\Gamma \varphi^2(\bm x)  \, d \bm x,
  \quad \|\mu\|_S^2 = \int_S \mu^2(\bm y)  \, d \bm y .
\]

\subsection{Class of nonnegative densities}

The Tikhonov regularization method \eqref{3.4} yields an approximate solution in the class of bounded solutions. 
We will work under the assumption that the density of the single-layer potential \eqref{3.1} is nonnegative.

We consider the inverse potential problem \eqref{2.1}, \eqref{2.2} under the a priori assumption $\rho(\bm{y}) \geq 0, \ \bm y \in D$.

Let us formulate an auxiliary boundary value problem in the domain $\Omega$ that includes the anomalies (Fig.~\ref{f-2}) ($D \subset \Omega$).
The function $w(\bm x), \bm x \in \Omega$ satisfies, like the volume potential $u(\bm x)$, the Poisson equation:

\begin{equation}\label{3.5}
  \Delta v = - \varrho(\bm x),
  \quad \bm x \in \Omega . 
\end{equation} 
On the boundary of the domain $\Omega$ we set the homogeneous Dirichlet condition:
\begin{equation}\label{3.6}
  w(\bm x) = 0, 
  \quad \bm x \in S .
\end{equation} 

For the solutions of the boundary value problem \eqref{3.5}, \eqref{3.6}, applying the third Green’s formula to points outside the extended domain $\Omega$ gives
\[
 \int_{D} \varrho(\bm y) G(\bm x,\bm y) d \bm y = - \int_S \frac{\partial v}{\partial n} (\bm y)  G(\bm x,\bm y) d \bm y ,
 \quad \bm x \in \Gamma , 
\] 
where $\bm n$ is the outward normal to $S$.
Thus, we have the representation of the volume potential \eqref{2.1} through the solution of the auxiliary problem \eqref{3.5}, \eqref{3.6}:
\begin{equation}\label{3.7}
 u(\bm x) = - \int_{S} \frac{\partial w}{\partial n} (\bm y)  G(\bm x,\bm y) d \bm y ,
 \quad \bm x \in \Gamma . 
\end{equation} 
Such a transformation to a problem of lower dimension is useful not only for solving the direct potential problem, but also for solving the ill-posed continuation problem toward the sources.

We will consider the class of inverse potential problems with sign-constant density.
For definiteness, we assume that $\varrho(\bm x) \geq  0, \ \bm x \in D$. 
The maximum principle for the boundary value problem \eqref{3.5}, \eqref{3.6} yields $w(\bm x) \geq  0, \ \bm x \in \Omega$ and therefore
\[
 \frac{\partial w}{\partial n} (\bm x) \geq 0, 
 \quad \bm x \in  \partial \Omega .
\] 
Taking into account \eqref{3.7}, we have a representation of the volume potential through a single-layer potential with nonnegative density. Thus, if $\rho(\bm{y}) \geq 0, \ \bm y \in D$, then $\mu(\bm{y}) \geq 0, \ \bm y \in S$.

These constraints are used in the approximate solution of the integral equation \eqref{3.3}. In this case,
\begin{equation}\label{3.8}
  \|\mathcal{A} \mu - \varphi \|_\Gamma^2 \rightarrow \min_{\mu \in K} ,
  \quad K = \{ s(\bm x) \ | \ s(\bm x) \geq 0, \ \bm x \in  S \} .
\end{equation}

\subsection{Discrete problem}

\begin{figure}[ht]
\centering
\includegraphics[width=0.5\textwidth]{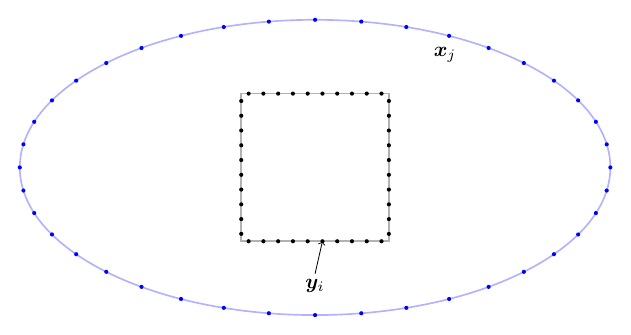}
\caption{Discretization.}
\label{f-3}
\end{figure}

The boundary $S$ is divided into small segments $S = \cup_{j=1}^N S_j$.
For the centers of the segments and their lengths we use the notation
\[
  \bm y_j, |S_j|, 
  \quad j = 1,2, \ldots, N .
\]
The potential is known at individual points on the boundary $\Gamma$:
\[
  \bm x_i, \quad i = 1,2, \ldots, M .
\]
The notation used is illustrated in Fig.~\ref{f-3}.

Quadrature approximation of the integral equation \eqref{3.3} leads to a system of linear algebraic equations
\begin{equation}\label{3.9}
  A v = f,
\end{equation}
where
\[
  v_j \approx \mu(\bm y_j),  \quad j = 1,2, \ldots, N ,
  \quad f_i = \varphi(\bm{x}_i) , \quad i = 1,2, \ldots, M .
\]
For the elements of the rectangular matrix $A$ we have
\[
 a_{ij} =   |S_j| G(\bm x_i, \bm y_j),
 \quad j = 1,2, \ldots, N , \quad i = 1,2, \ldots, M .
\]

In the simplest case, the least squares method is used for the approximate solution of \eqref{3.9}.
In these conditions, the approximate solution is defined from 
\begin{equation}\label{3.10}
  \|A v - f \|^2 \rightarrow \min_{v} .
\end{equation}

Accounting for measurement errors at the observation points $\bm x_i, \quad i = 1,2, \ldots, M$ leads to the need to solve equation \eqref{3.9} with an inexact right-hand side, where $\widetilde{f} \approx f$.
A stable approximate solution can be obtained by using the Tikhonov regularization method.
For the system of linear equations \eqref{3.9}, the approximate solution for a given regularization parameter $\alpha$ is defined from 
\begin{equation}\label{3.11}
  \|A v - \widetilde{f} \|^2 + \alpha \|v\|^2 \rightarrow \min_{v} .
\end{equation}
The choice of $\alpha$ is consistent with the error in the right-hand side, using, for example, the discrepancy principle.

We consider the problem of approximating the observed potential field by a single-layer potential in the class of nonnegative densities $\mu(\bm{y}) \geq 0, \ \bm y \in S$.
The corresponding discrete problem is formulated as follows:
\begin{equation}\label{3.12}
  \|A v - f \|^2 \rightarrow \min_{v > 0} .
\end{equation}
The computational implementation of \eqref{3.12} is based on the iterative NNLS (Nonnegative Least Squares) algorithm.

\section{Computational algorithms and numerical experiments}

The computational identification of the source region $\Omega$ is based on approximating the volumetric integral by a surface integral on the boundary $S$.
After discretization, the problem reduces to solving a system of linear algebraic equations, for which various regularization algorithms are applied.
Below, three approaches are considered: the least squares method, Tikhonov regularization, and the NNLS algorithm.

\subsection{Test problem}

Consider a model inverse potential problem with synthetic data.
The potential is measured on the boundary of an ellipse $\Gamma$ centered at $(0,0)$ with semi-axes 2 and 1:
\[
  \Gamma = \{\bm x \ | \ 0.25 (x^{(1)})^2 + (x^{(2)})^2 = 1 \}.
\]
The boundary is discretized uniformly along the angle into observation points $\bm x_i, \ i = 1,2,\dots,M$.

The potential source region $\Omega$ is a unit square.
Its boundary $S$ is divided into $N_1$ equal segments horizontally and $N_2$ vertically.
The data used in the inverse problem are computed analytically for the volumetric potential.
Assume $\rho(\bm y) = 1, \ \bm y \in D$, where $D$ consists of two disks: radius 0.1 centered at $(-0.2,0)$ and radius 0.05 centered at $(0.2,-0.2)$.

\begin{figure}[ht]
\centering
\includegraphics[width=0.75\textwidth]{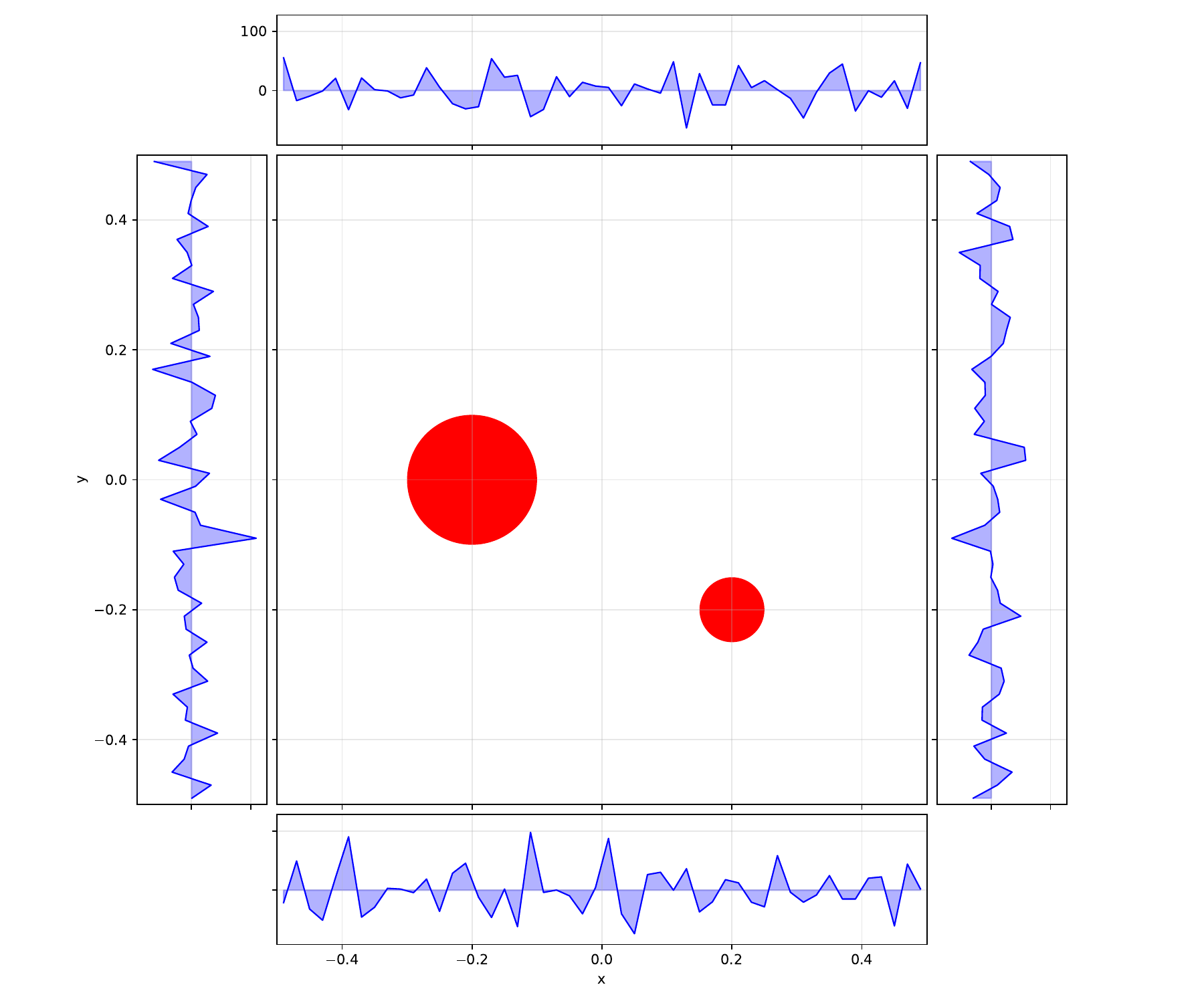}
\caption{Solution of the problem for $N_1 = N_2 = 50$.}
\label{f-4}
\end{figure}

\subsection{Least squares method}

\begin{figure}[ht]
\centering
\includegraphics[width=0.75\textwidth]{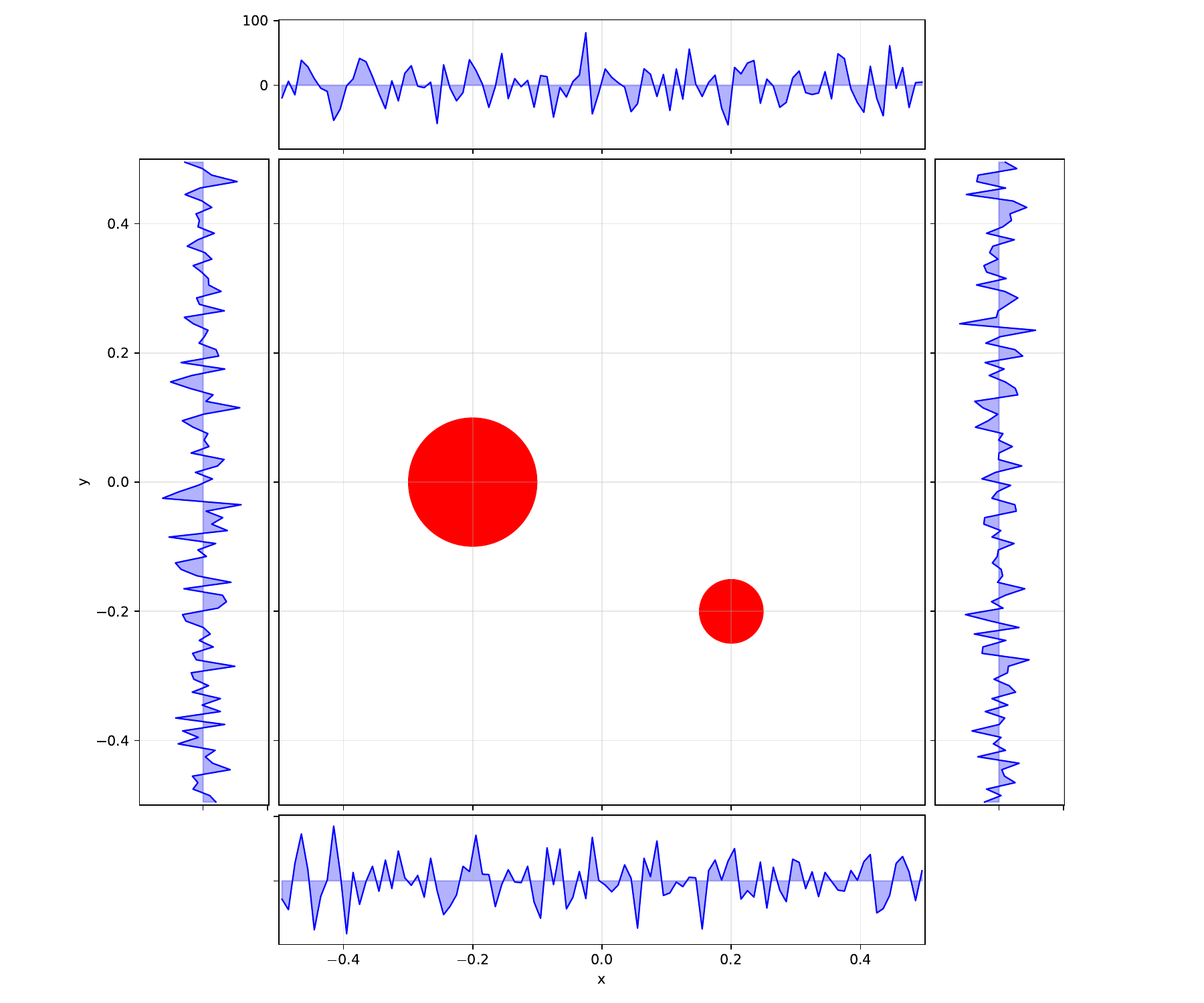}
\caption{Solution of the problem for $N_1 = N_2 = 100$.}
\label{f-5}
\end{figure}

\begin{figure}[ht]
\centering
\includegraphics[width=0.75\textwidth]{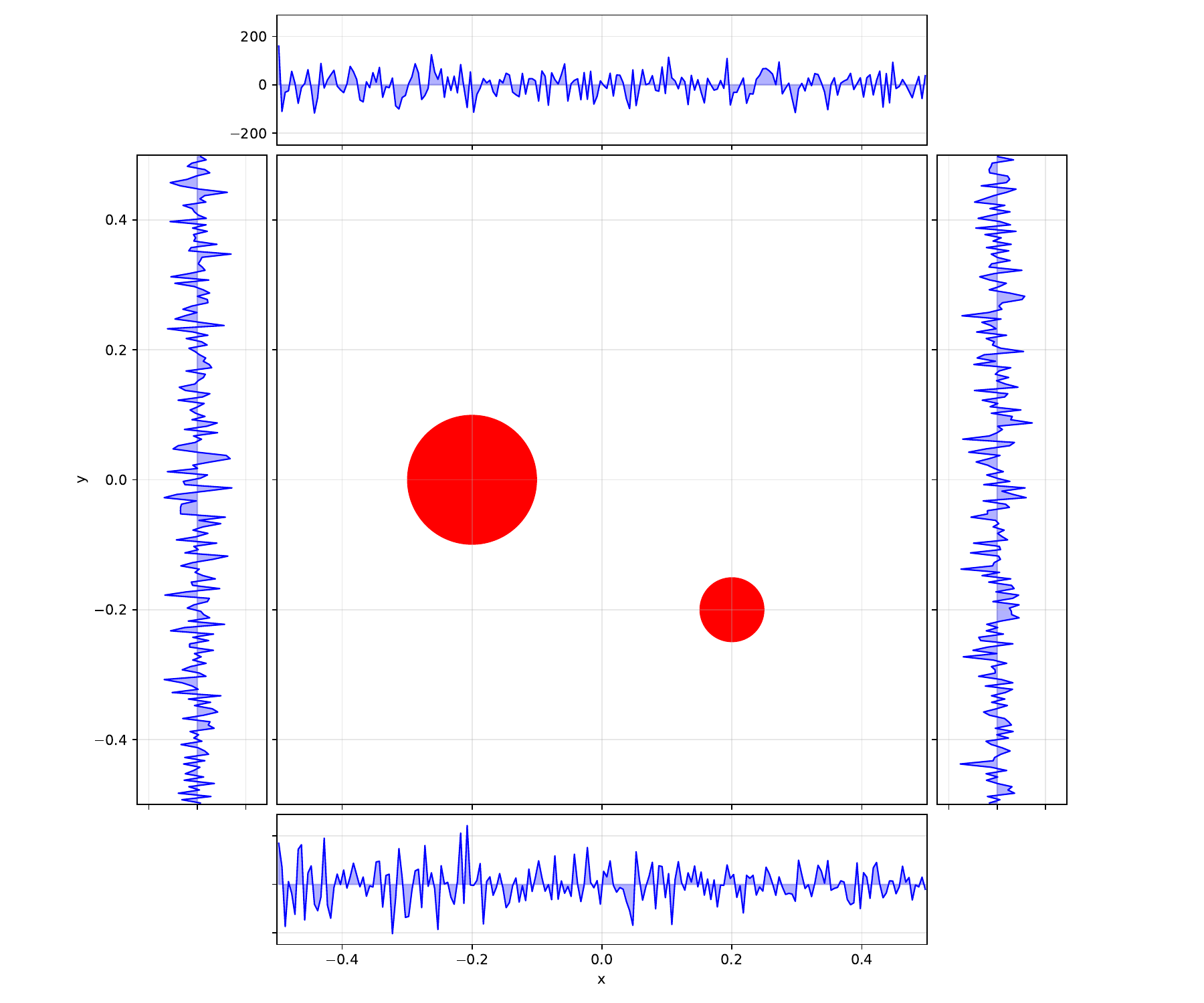}
\caption{Solution of the problem for $N_1 = N_2 = 200$.}
\label{f-6}
\end{figure}

After discretizing $S$ and $\Gamma$, the inverse problem is formulated as a rectangular linear system of size $N \times M$, $N = 2 (N_1 + N_2)$.
Its approximate solution is obtained using the standard least squares method \eqref{3.10}.
Density reconstruction results using synthetic data are shown in Figures~\ref{f-4}--\ref{f-6}.

The recovered simple-layer potential density exhibits sign changes.
The amplitude of oscillations depends strongly on the discretization of $S$, although the potential approximation on $\Gamma$ remains accurate.

\begin{figure}[ht]
\centering
\includegraphics[width=0.75\textwidth]{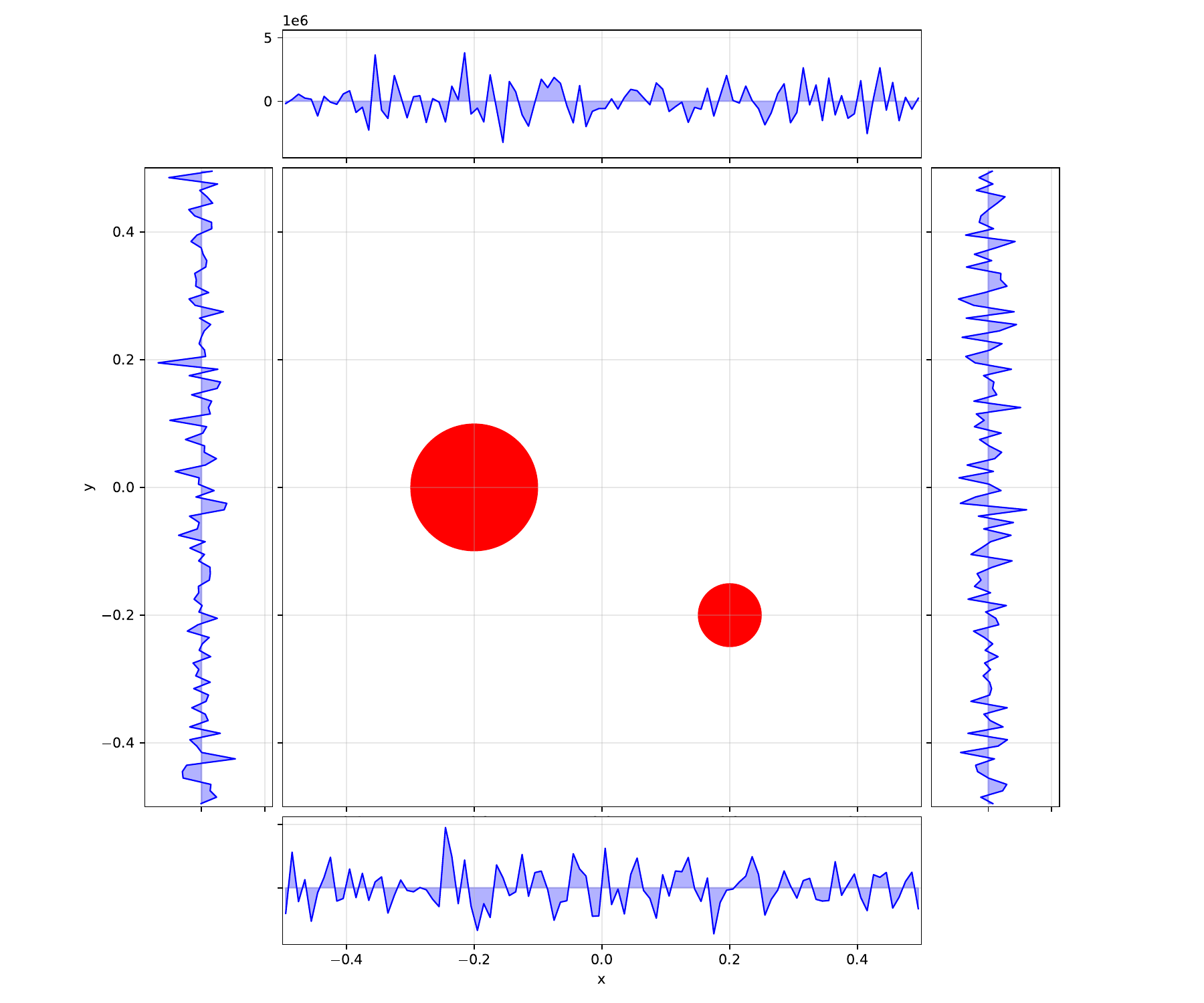}
\caption{Solution of the problem for $\delta = 0.05$.}
\label{f-7}
\end{figure}

\begin{figure}[ht]
\centering
\includegraphics[width=0.75\textwidth]{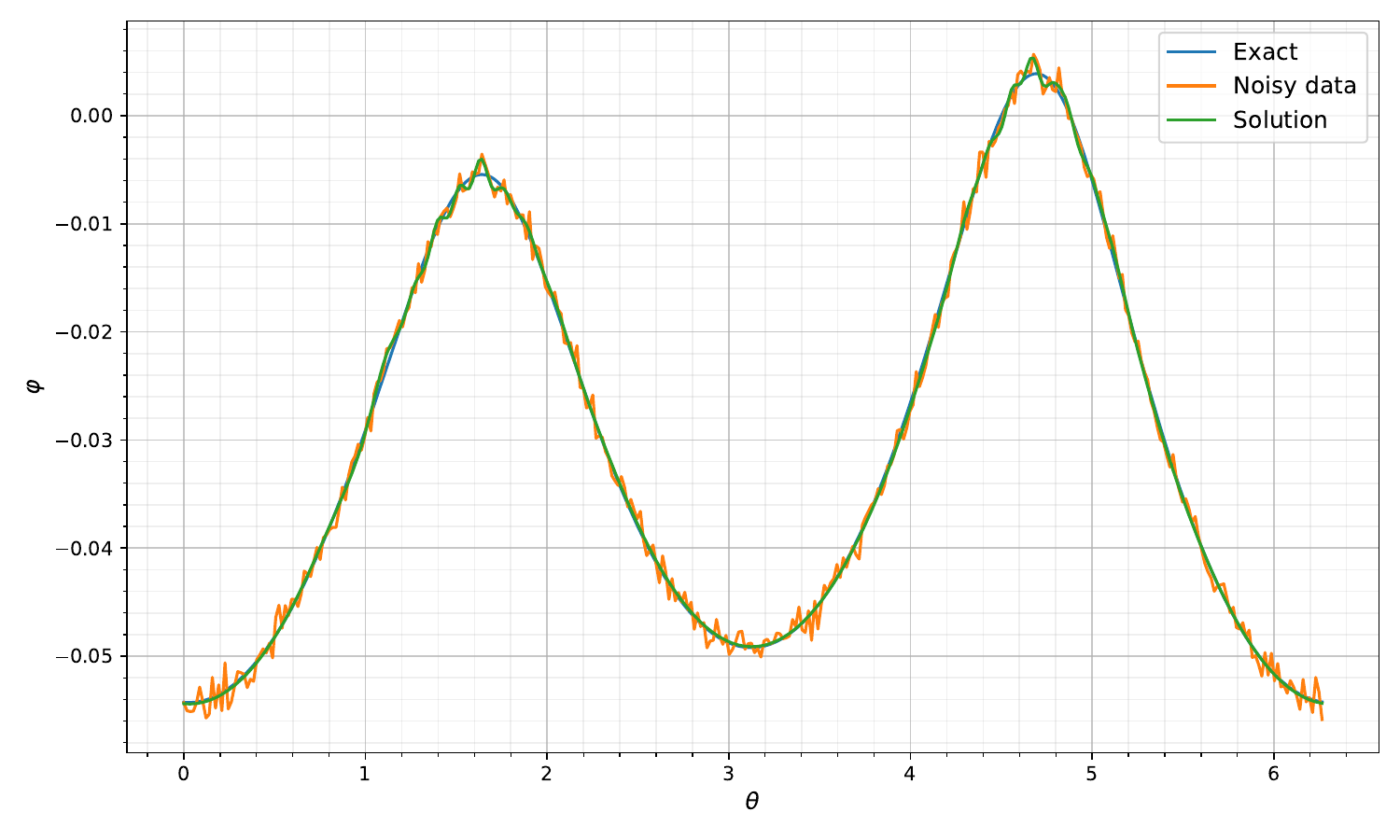}
\caption{Potential at observation points for $\delta = 0.05$.}
\label{f-8}
\end{figure}

With noisy data, defined as
\[
  \widetilde{f}_i = f_i + \delta s(f) \sigma_i,
\]
where $\delta$ is the noise amplitude, $s(f)$ is the standard deviation, and $\sigma_i \sim N(0,1)$,
the amplitude and frequency of oscillations increase (Figures~\ref{f-7}, \ref{f-8}).

\textbf{Conclusion.} The least squares method provides a good approximation on $\Gamma$, but produces strong oscillations and sign changes in the density, especially for noisy data.

\clearpage

\subsection{Tikhonov regularization}

\begin{figure}[ht]
\centering
\includegraphics[width=0.75\textwidth]{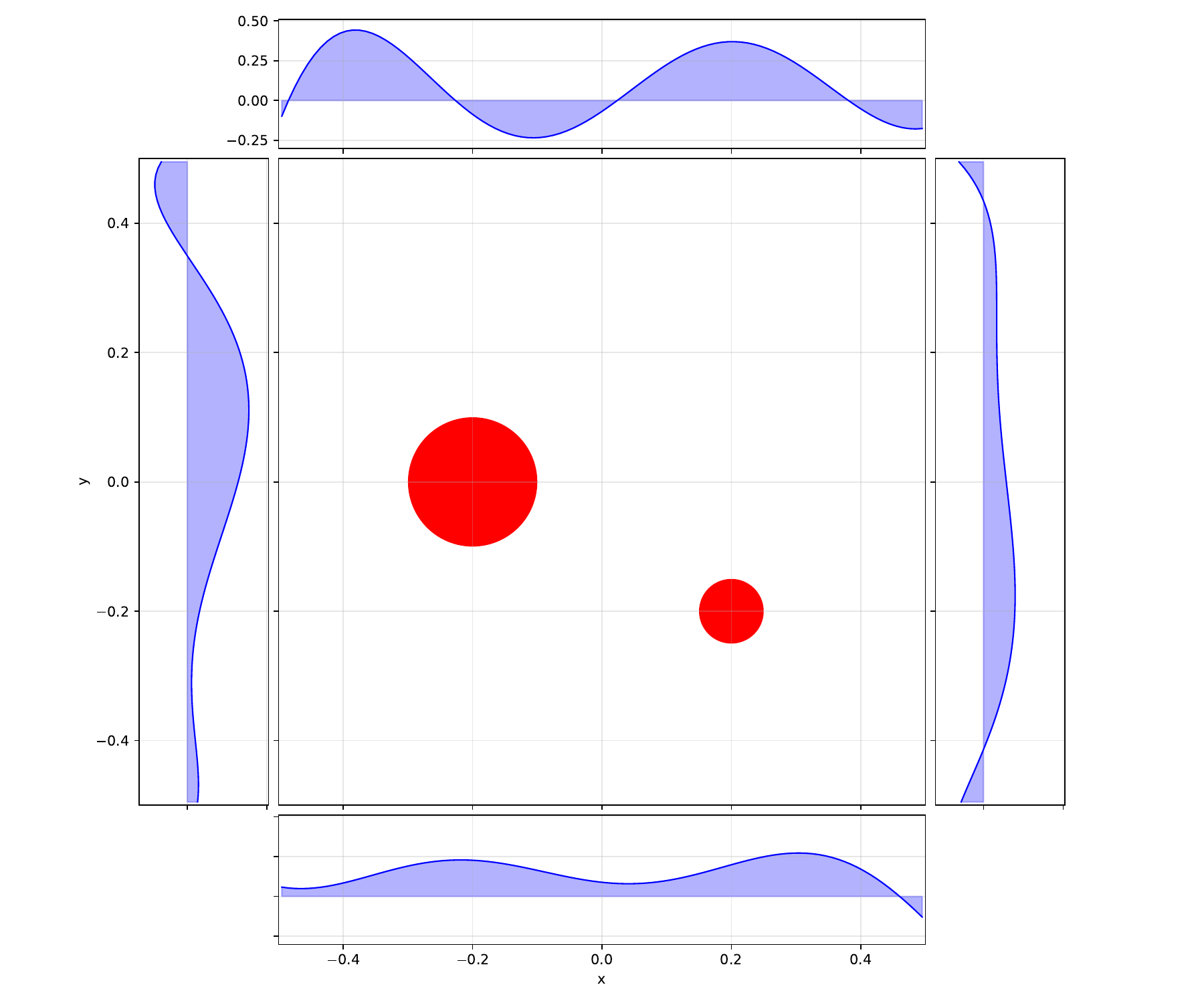}
\caption{Solution of the problem for $\alpha = 10^{-7}$ (noise level $\delta = 0.05$).}
\label{f-9}
\end{figure}

\begin{figure}[ht]
\centering
\includegraphics[width=0.75\textwidth]{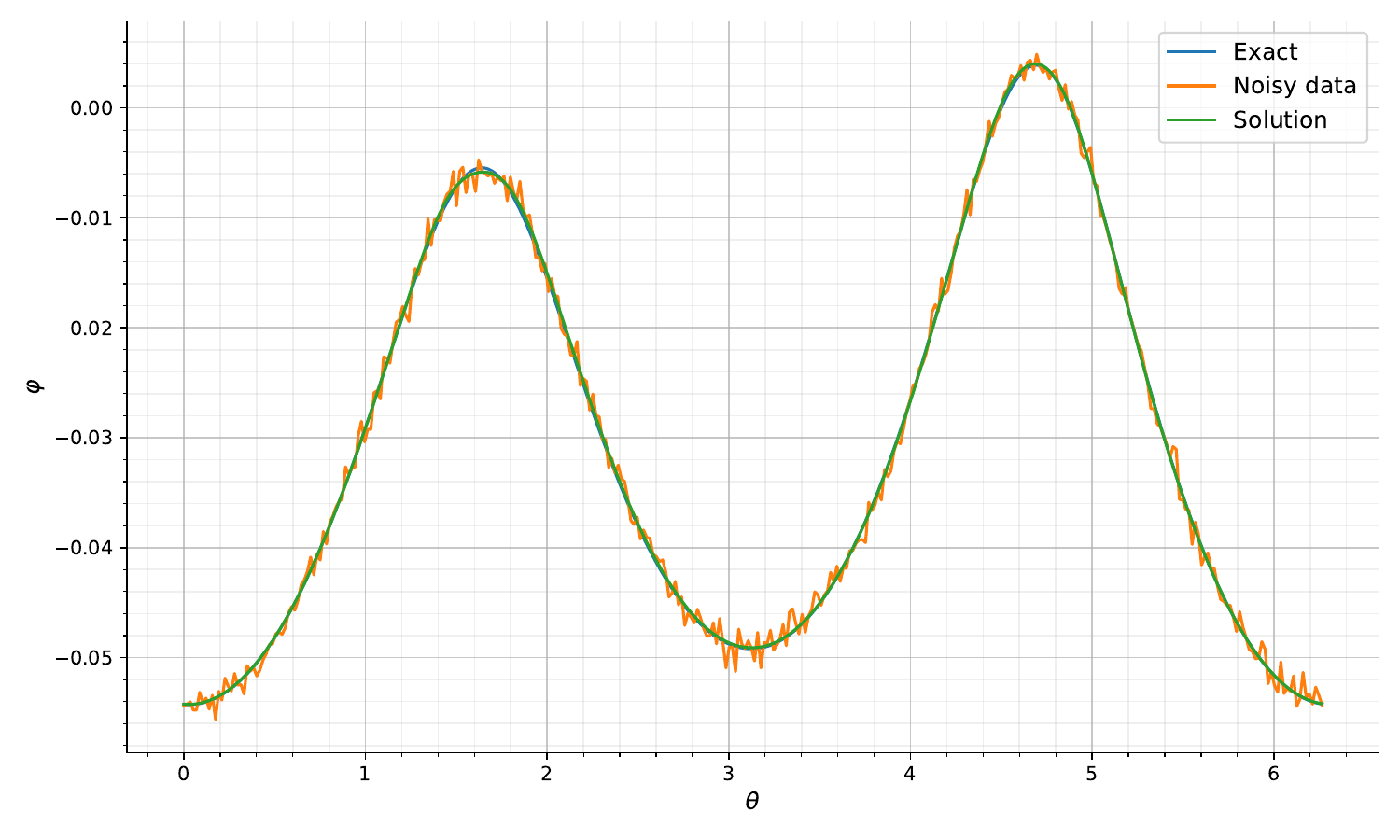}
\caption{Potentials at observation points for $\alpha = 10^{-7}$ (noise level $\delta = 0.05$).}
\label{f-10}
\end{figure}

\begin{figure}[ht]
\centering
\includegraphics[width=0.75\textwidth]{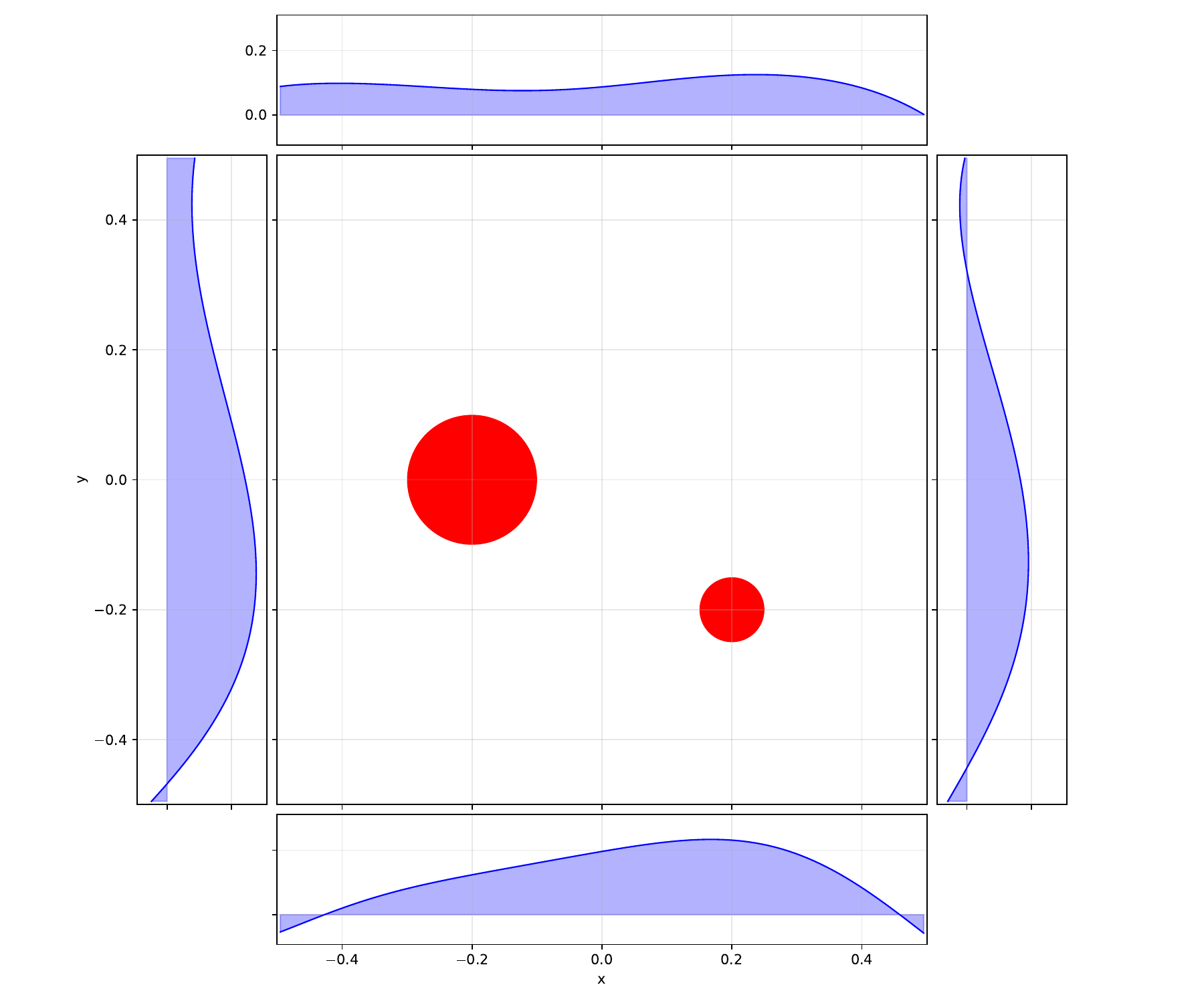}
\caption{Solution of the problem for $\alpha = 10^{-6}$ (noise level $\delta = 0.05$).}
\label{f-11}
\end{figure}

\begin{figure}[ht]
\centering
\includegraphics[width=0.75\textwidth]{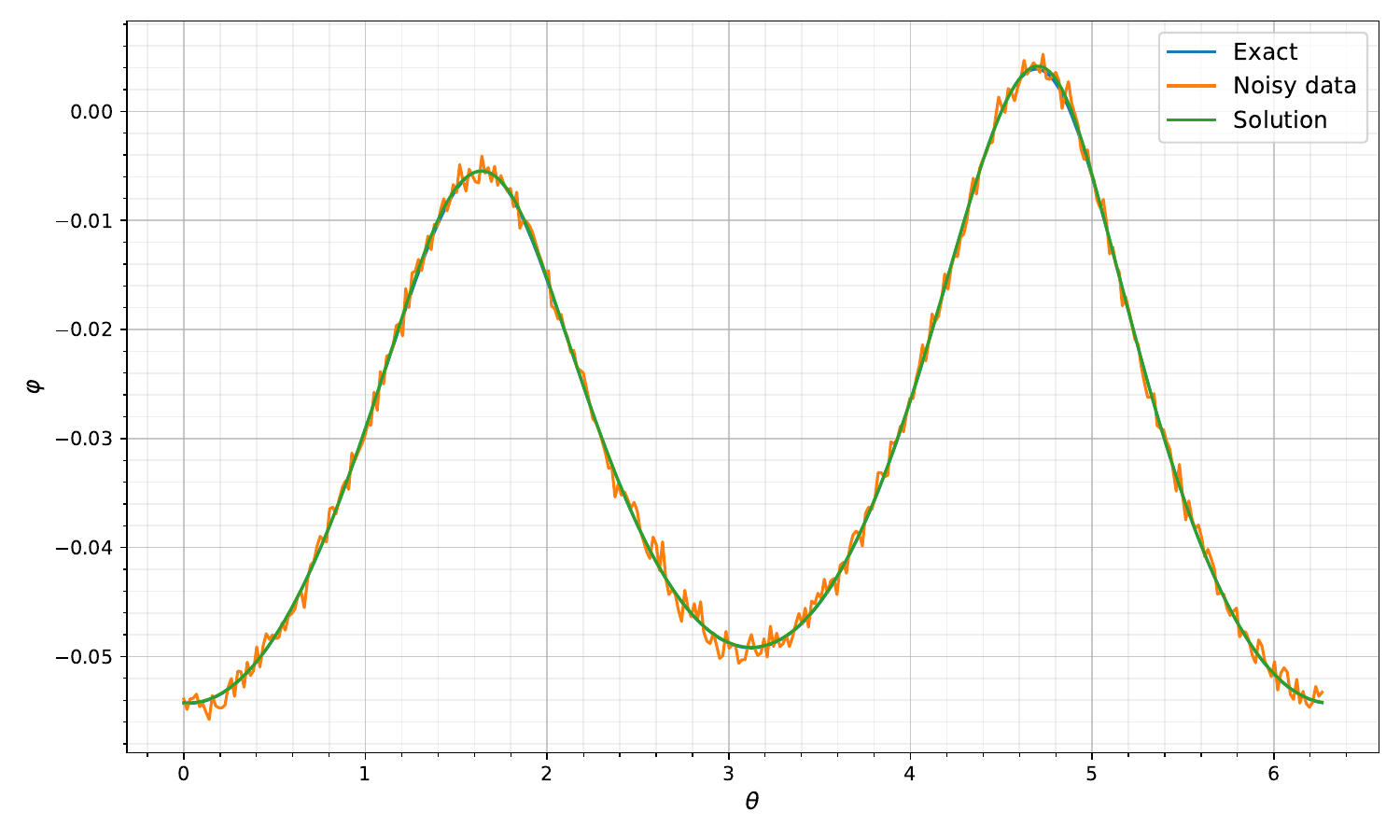}
\caption{Potentials at observation points for $\alpha = 10^{-6}$ (noise level $\delta = 0.05$).}
\label{f-12}
\end{figure}

Using Tikhonov regularization, the approximate solution is obtained from \eqref{3.11}.
An example with noise level $\delta = 0.05$ and regularization parameter $\alpha = 10^{-7}$ is shown in Figures~\ref{f-9}, \ref{f-10}.
Oscillations are significantly suppressed, highlighting the minimum-norm solution.
Increasing $\alpha$ results in further smoothing (Figures~\ref{f-11}, \ref{f-12}).

\textbf{Conclusion.} Tikhonov regularization stabilizes the solution and reduces oscillations, but smoothing may diminish details in the density distribution.

\clearpage

\subsection{NNLS algorithm}

\begin{figure}[ht]
\centering
\includegraphics[width=0.75\textwidth]{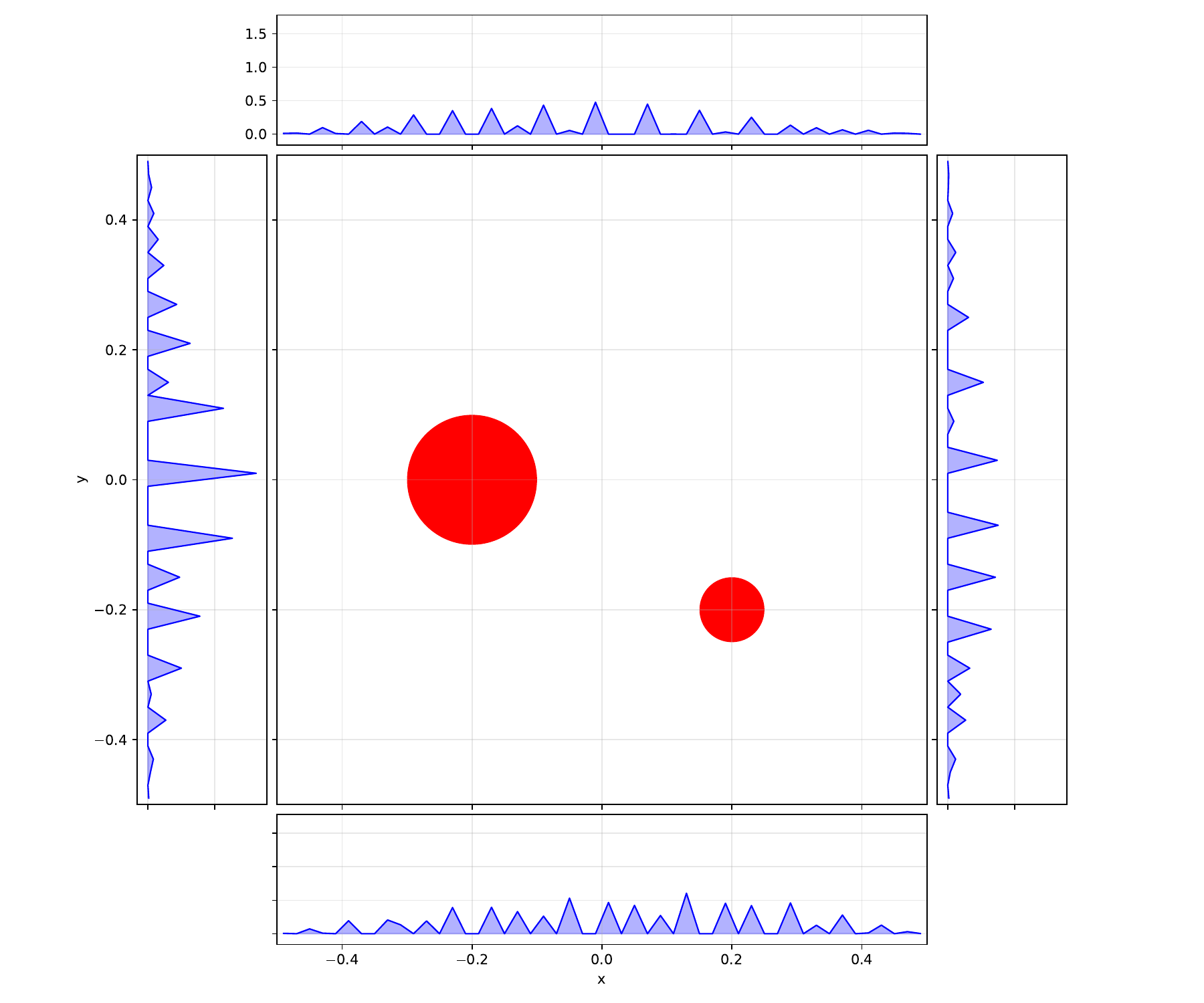}
\caption{Solution of the problem using the NNLS algorithm for $N_1 = N_2 = 50$.}
\label{f-13}
\end{figure}

\begin{figure}[ht]
\centering
\includegraphics[width=0.75\textwidth]{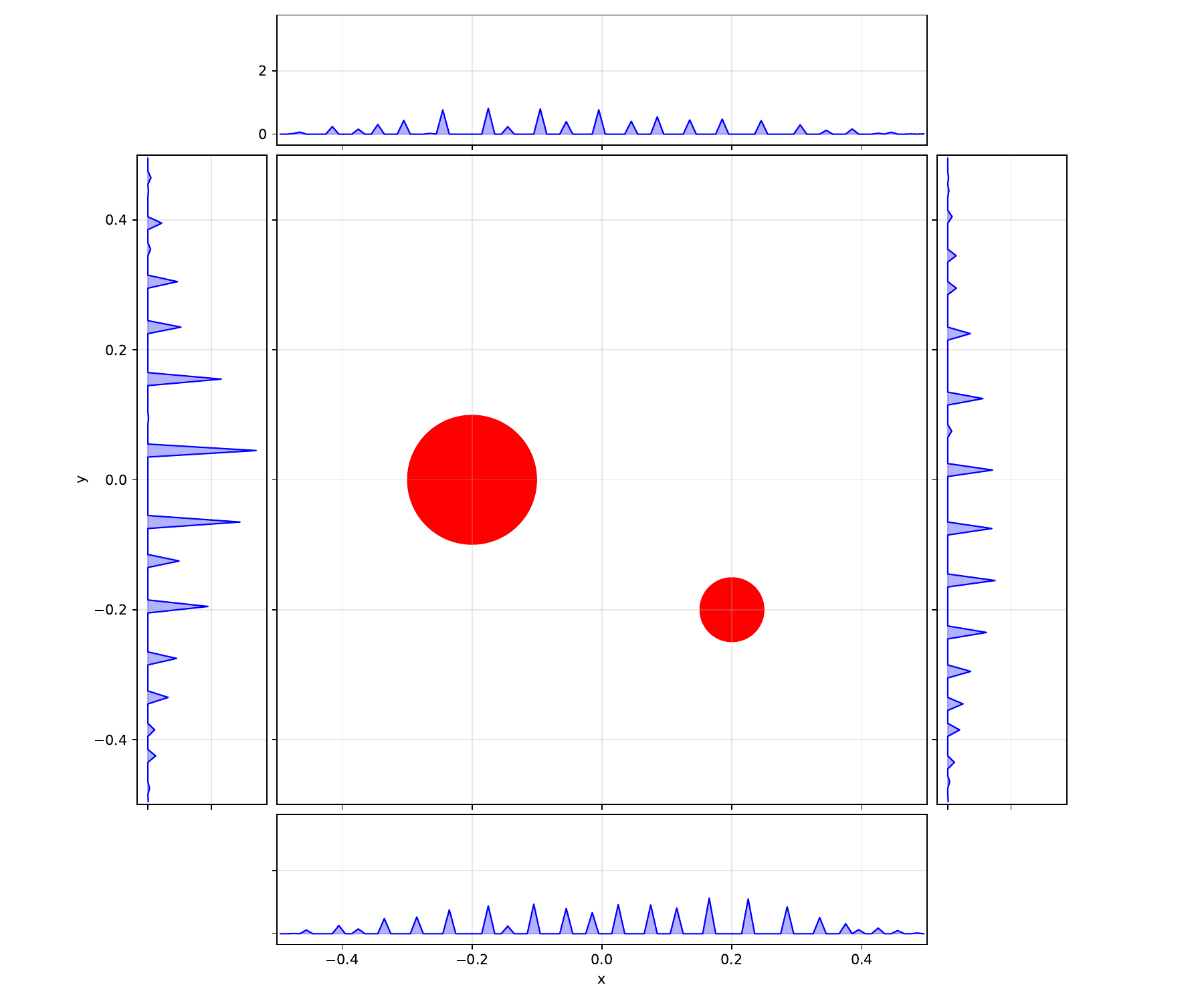}
\caption{Solution of the problem using the NNLS algorithm for $N_1 = N_2 = 100$.}
\label{f-14}
\end{figure}

\begin{figure}[ht]
\centering
\includegraphics[width=0.75\textwidth]{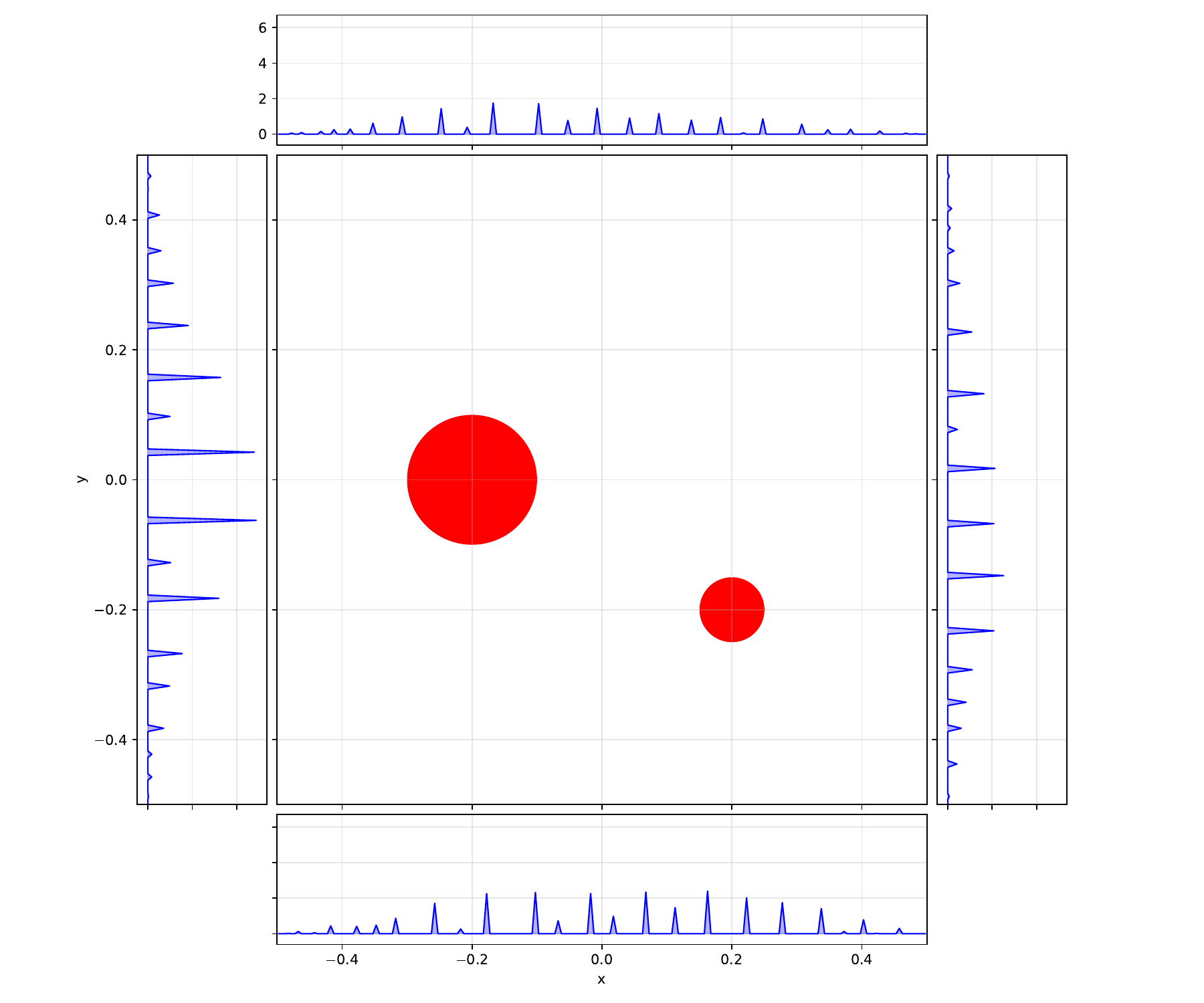}
\caption{Solution of the problem using the NNLS algorithm for $N_1 = N_2 = 200$.}
\label{f-15}
\end{figure}

Approximating the observed potential by a simple-layer potential is considered in the class of nonnegative densities $\mu(\bm y) \geq 0, \ \bm y \in S$.
The discrete problem is formulated as
\[
  \|A v - f \|^2 \rightarrow \min_{v > 0}.
\]
It is solved using the Python implementation of NNLS (function \texttt{nnls()} from the SciPy library).
Density reconstruction for different discretizations is shown in Figures~\ref{f-13}--\ref{f-15}.

\begin{figure}[ht]
\centering
\includegraphics[width=0.75\textwidth]{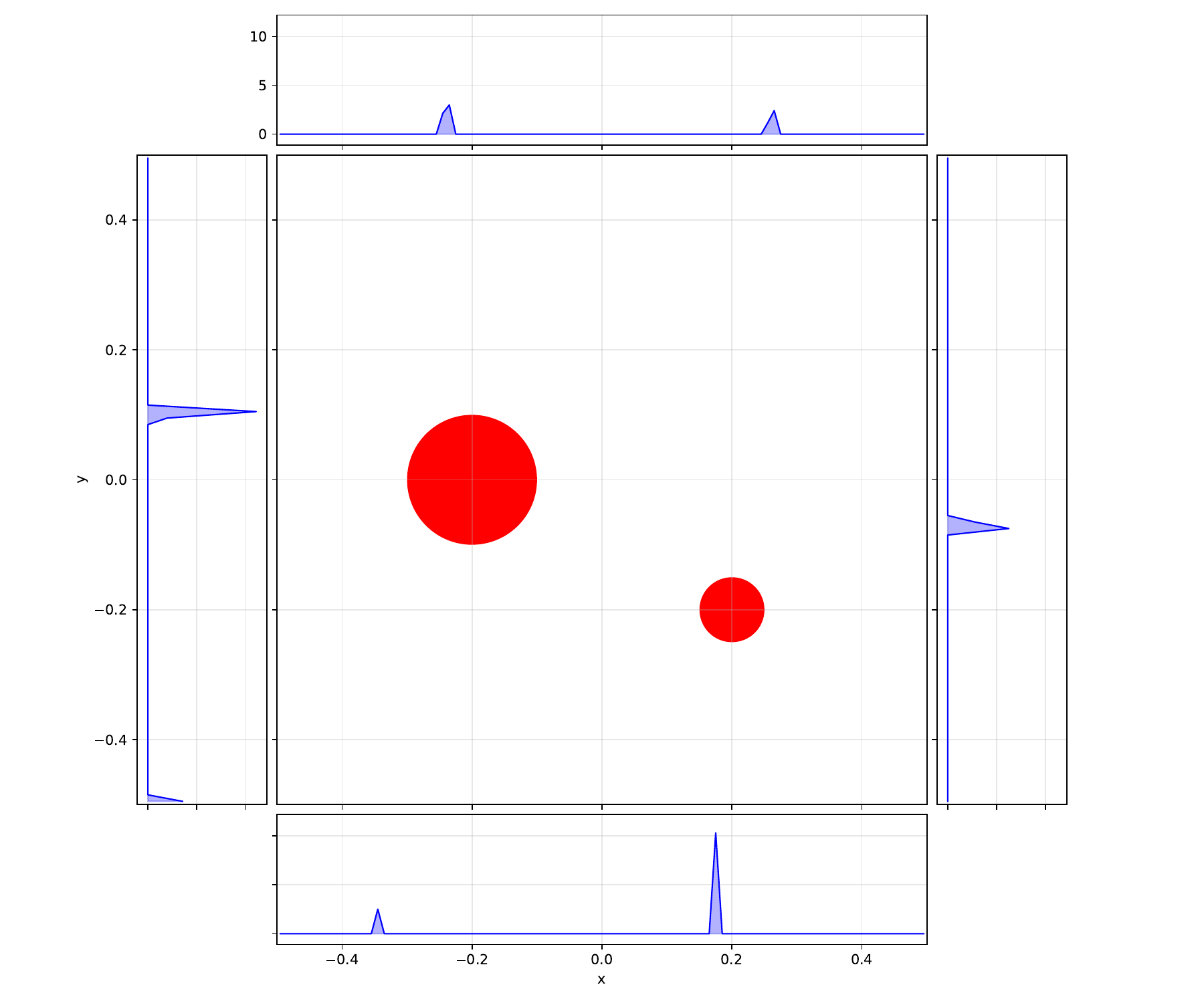}
\caption{Solution of the problem for $\delta = 0.05$ using the NNLS algorithm.}
\label{f-16}
\end{figure}

\begin{figure}[ht]
\centering
\includegraphics[width=0.75\textwidth]{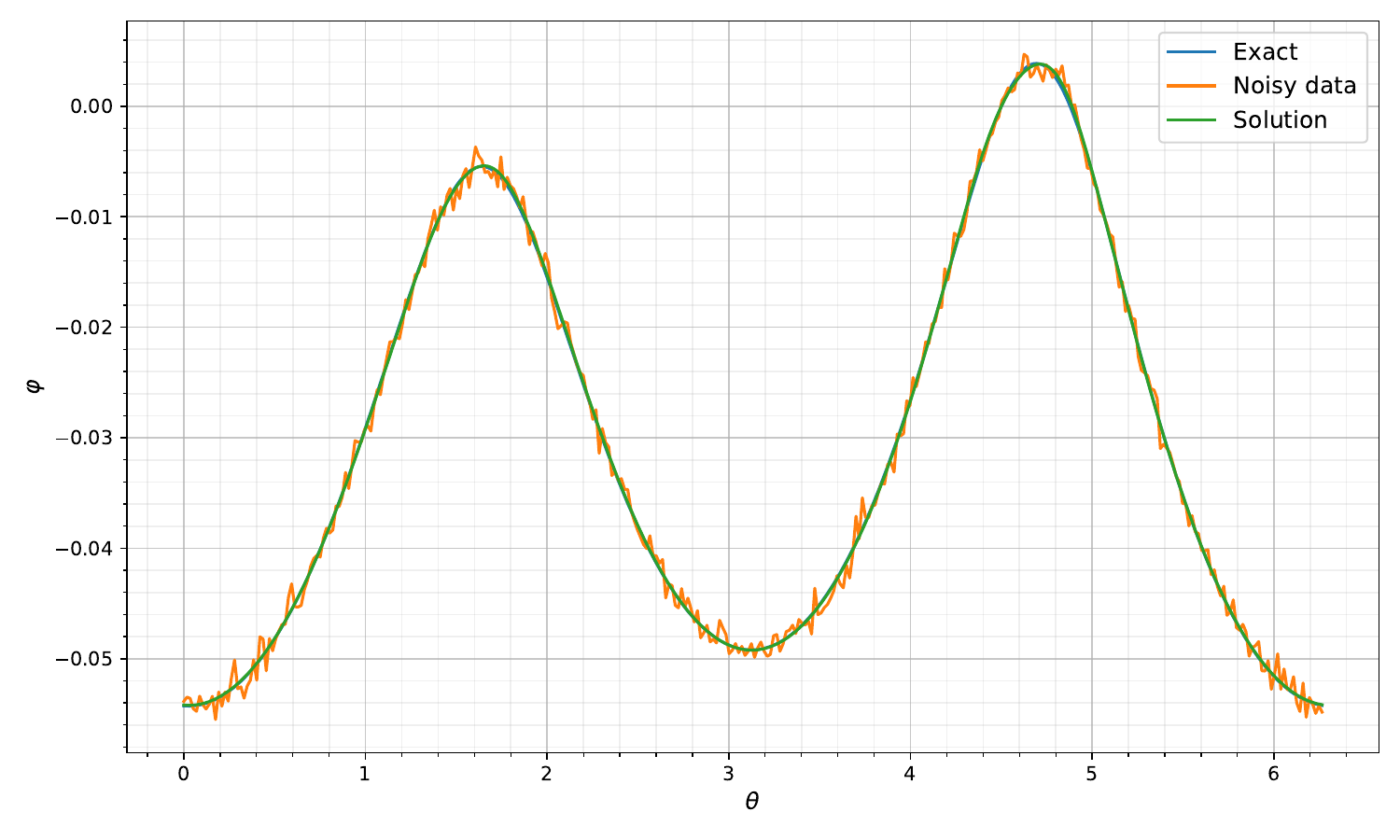}
\caption{Potentials at observation points for $\delta = 0.05$ using the NNLS algorithm.}
\label{f-17}
\end{figure}

For noisy data ($\delta = 0.05$), the approximation remains accurate (Figures~\ref{f-16}, \ref{f-17}).
Quality is achieved through fewer nonzero components in $v$, corresponding to localization of sources along specific segments of $S$.

\textbf{Conclusion.} NNLS effectively exploits the nonnegativity prior, eliminates oscillations, and improves robustness under noisy conditions.

\clearpage

\subsection{Data for $D \not\subset \Omega$}

The source region $D$ cannot be uniquely determined due to the fundamental ill-posedness of the inverse problem.
Hence, the exact distribution of the equivalent simple-layer potential density is less important.
The critical question is whether the computational data allow us to determine if $D$ lies entirely within a given rectangle $\Omega$, and whether $D$ can be localized by varying $\Omega$.

\begin{figure}[ht]
\centering
\includegraphics[width=1\textwidth]{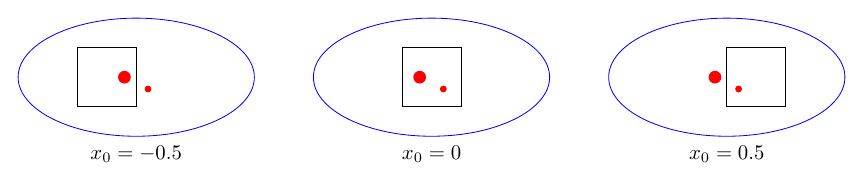}
\caption{Computational scheme for different positions of the rectangle $S$.}
\label{f-18}
\end{figure}

\begin{figure}[ht]
\centering
\includegraphics[width=0.75\textwidth]{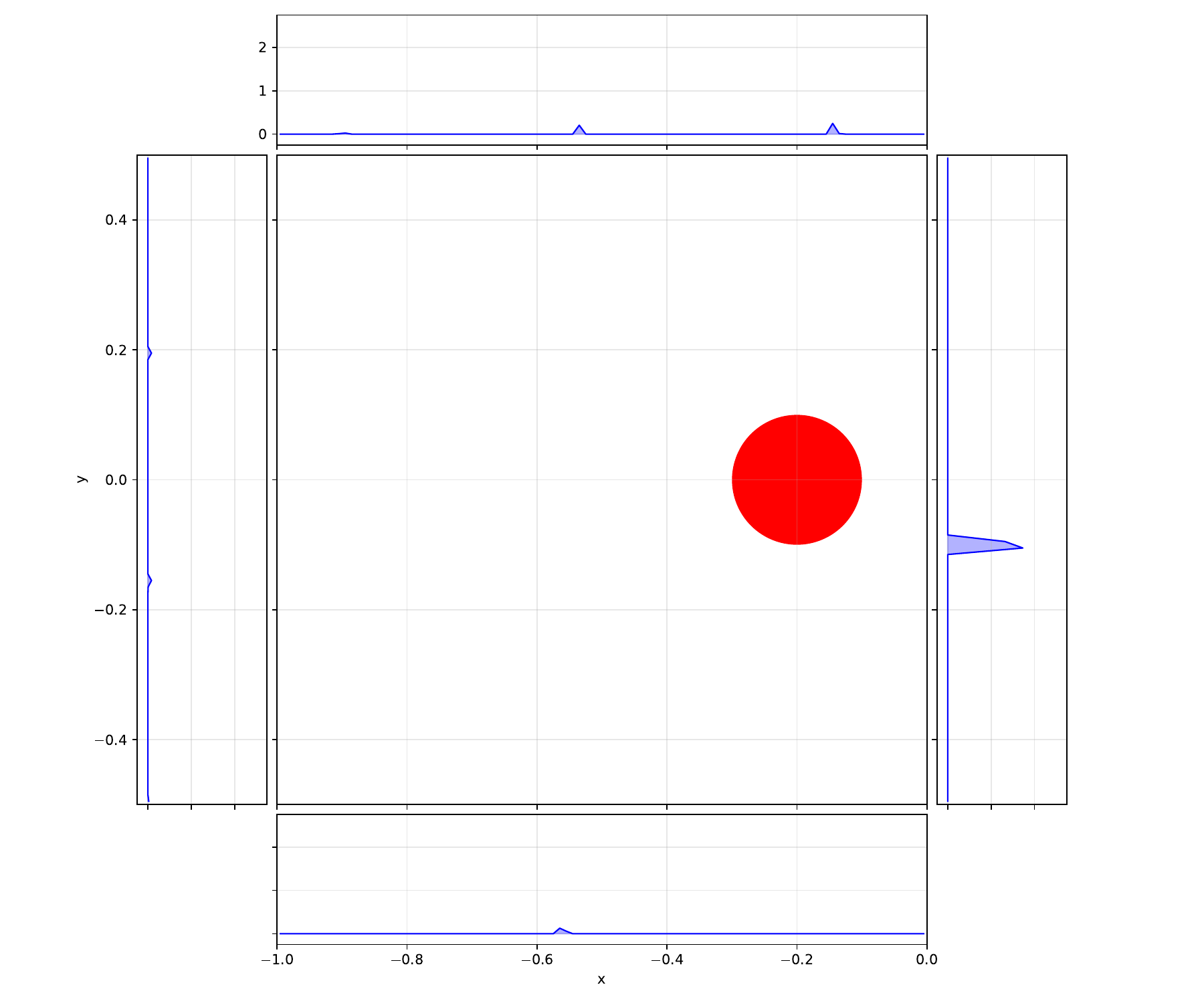}
\caption{Solution of the problem for $D_2 \notin \Omega$ ($x_0 = -0.5$).}
\label{f-19}
\end{figure}

\begin{figure}[ht]
\centering
\includegraphics[width=0.75\textwidth]{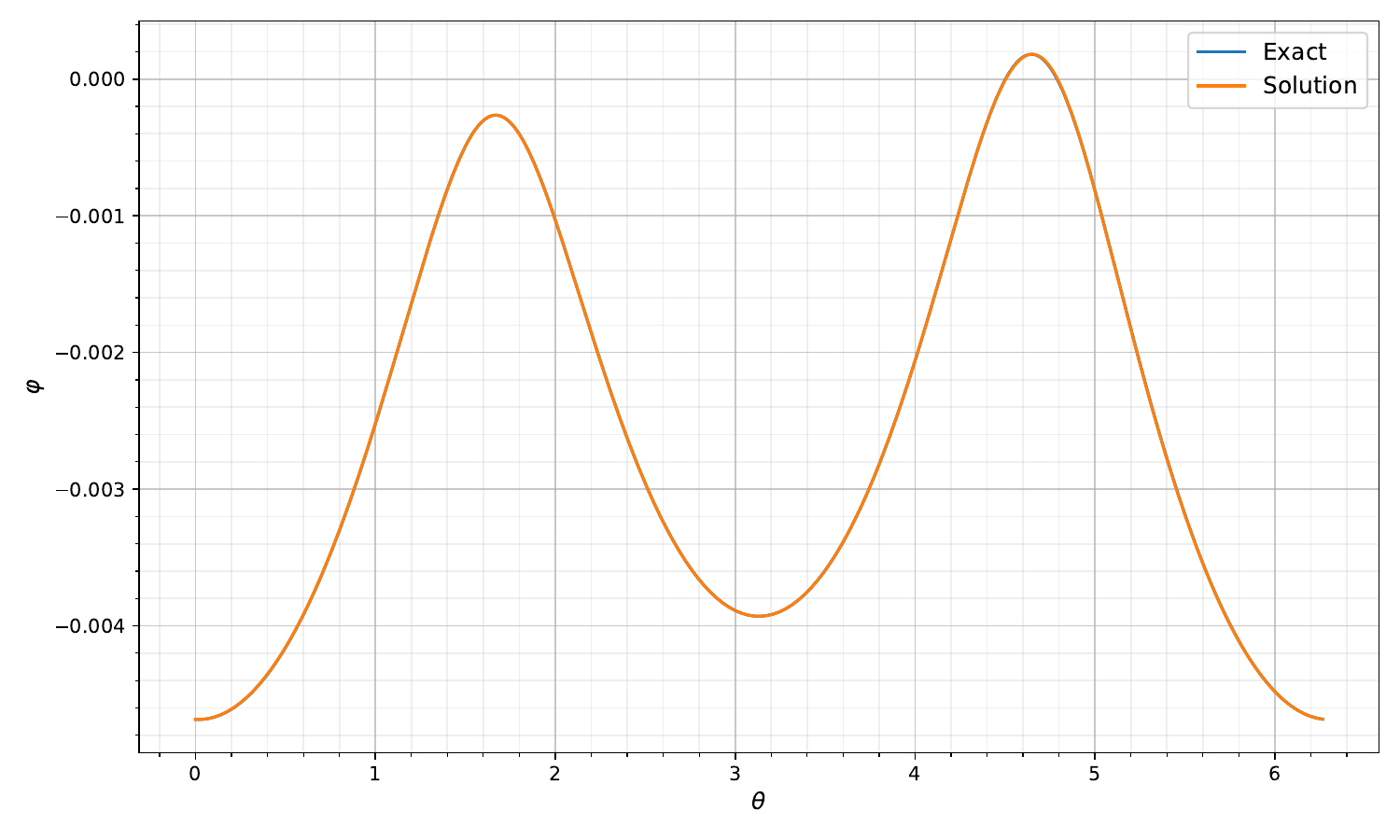}
\caption{Potentials at observation points for $D_2 \notin \Omega$ ($x_0 = -0.5$).}
\label{f-20}
\end{figure}

\begin{figure}[ht]
\centering
\includegraphics[width=0.75\textwidth]{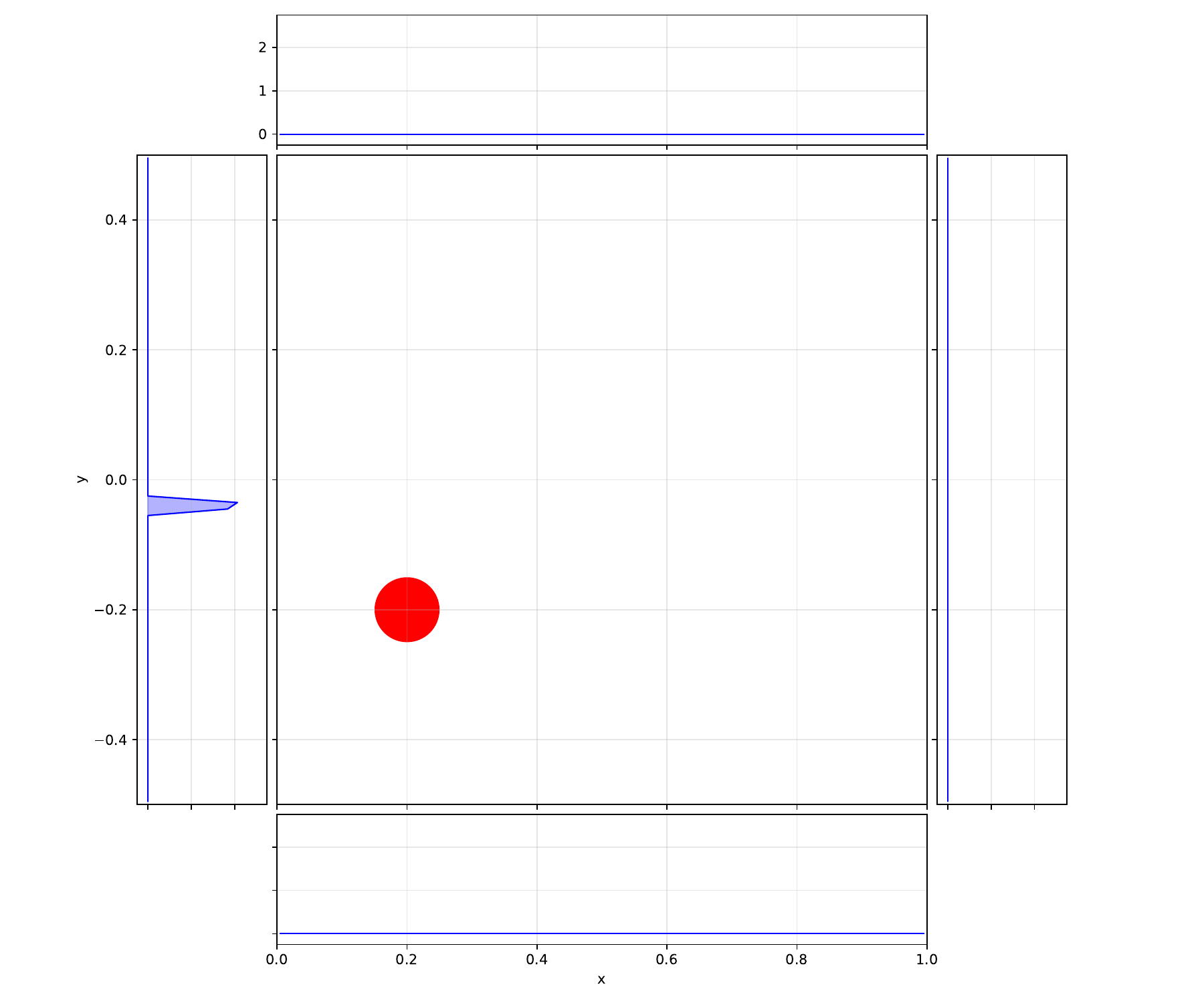}
\caption{Solution of the problem for $D_1 \notin \Omega$ ($x_0 = 0.5$).}
\label{f-21}
\end{figure}

\begin{figure}[ht]
\centering
\includegraphics[width=0.75\textwidth]{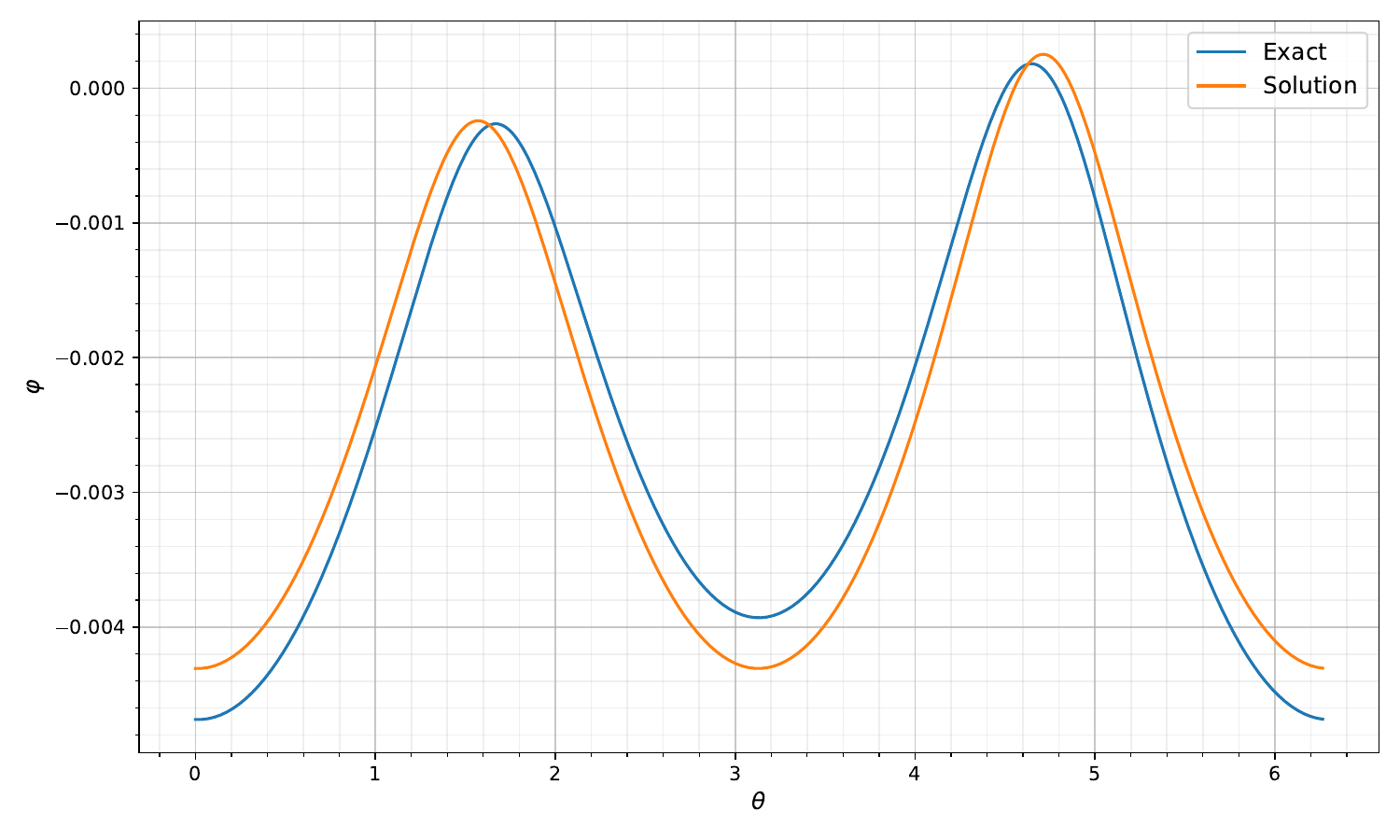}
\caption{Potentials at observation points for $D_1 \notin \Omega$ ($x_0 = 0.5$).}
\label{f-22}
\end{figure}

Calculations are performed for nonnegative densities when $D$ is not entirely inside the rectangle $\Omega$, moving the rectangle center along $x_0$ (Figure~\ref{f-18}).
Figures~\ref{f-19}--\ref{f-22} illustrate the effect of parts of $D$ lying outside $\Omega$.
The maximum simple-layer density shifts to the boundary, and the residual increases, indicating loss of approximation accuracy on $\Gamma$.
This provides a practical indication of the source region location.

\clearpage

\subsection{NNDD algorithm}

The Nonnegative Density Domain (NNDD) algorithm for source localization proceeds as follows:
\begin{itemize}
\item Specify the size of the window $S$ for the expected source region.
\item Vary the window position (e.g., rectangle center $x_0, y_0$) and compute the fit using NNLS for each position.
\item Identify the source region by the minimum residual.
\end{itemize}

\begin{figure}[ht]
\centering
\includegraphics[width=0.75\textwidth]{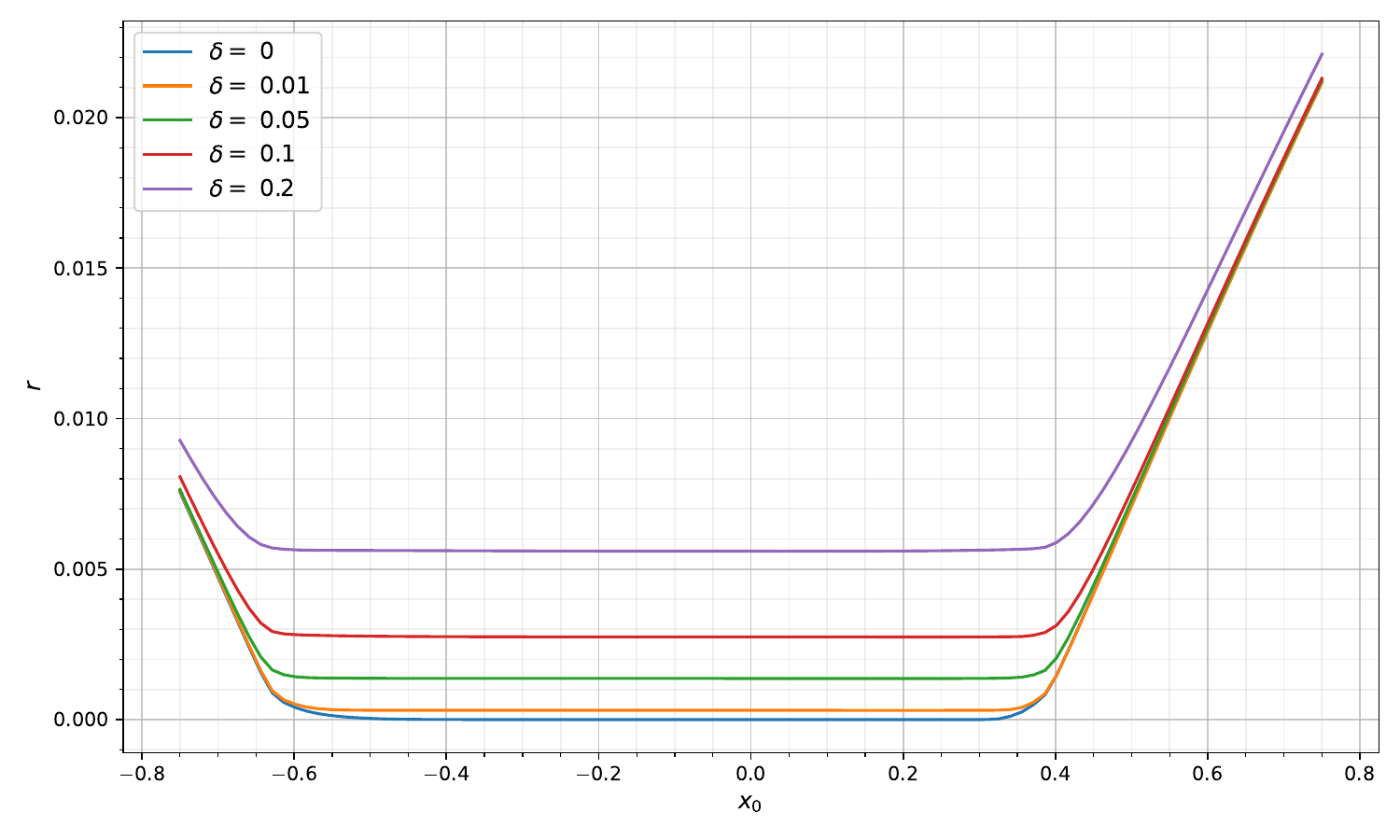}
\caption{Residuals for different positions of $S$.}
\label{f-23}
\end{figure}

\begin{figure}[ht]
\centering
\includegraphics[width=0.75\textwidth]{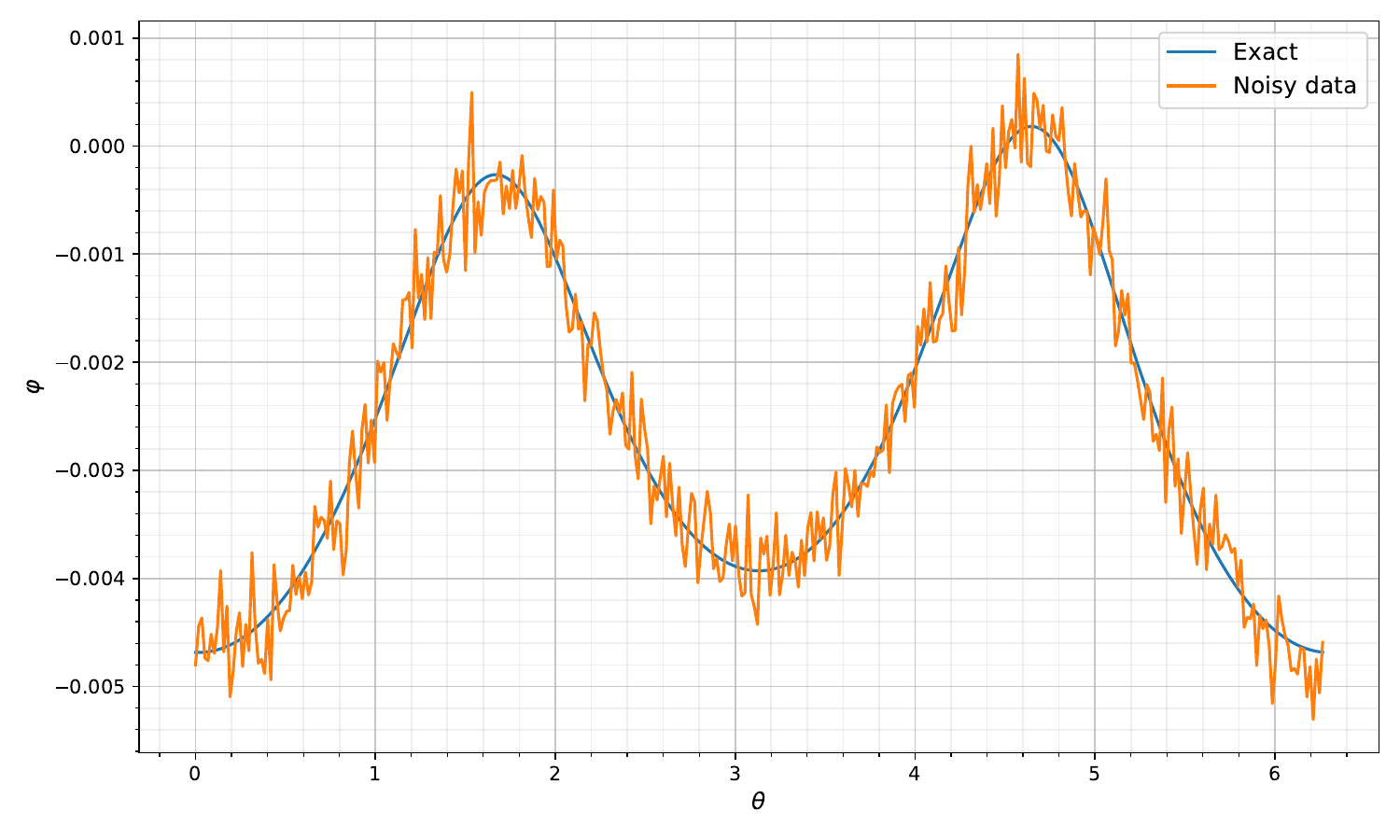}
\caption{Data perturbation for $\delta = 0.2$.}
\label{f-24}
\end{figure}

For the test problem, the rectangle center is varied along $x$ ($y_0 = 0$).
Figure~\ref{f-23} shows residual variation with $x_0$ for several noise levels.
A smooth residual plateau occurs when sources are inside the window $S$, allowing geometric localization of heterogeneities.
Even at high noise levels ($\delta = 0.2$, Figure~\ref{f-24}), the localization region remains reasonably robust.

\begin{figure}[ht]
\centering
\includegraphics[width=0.75\textwidth]{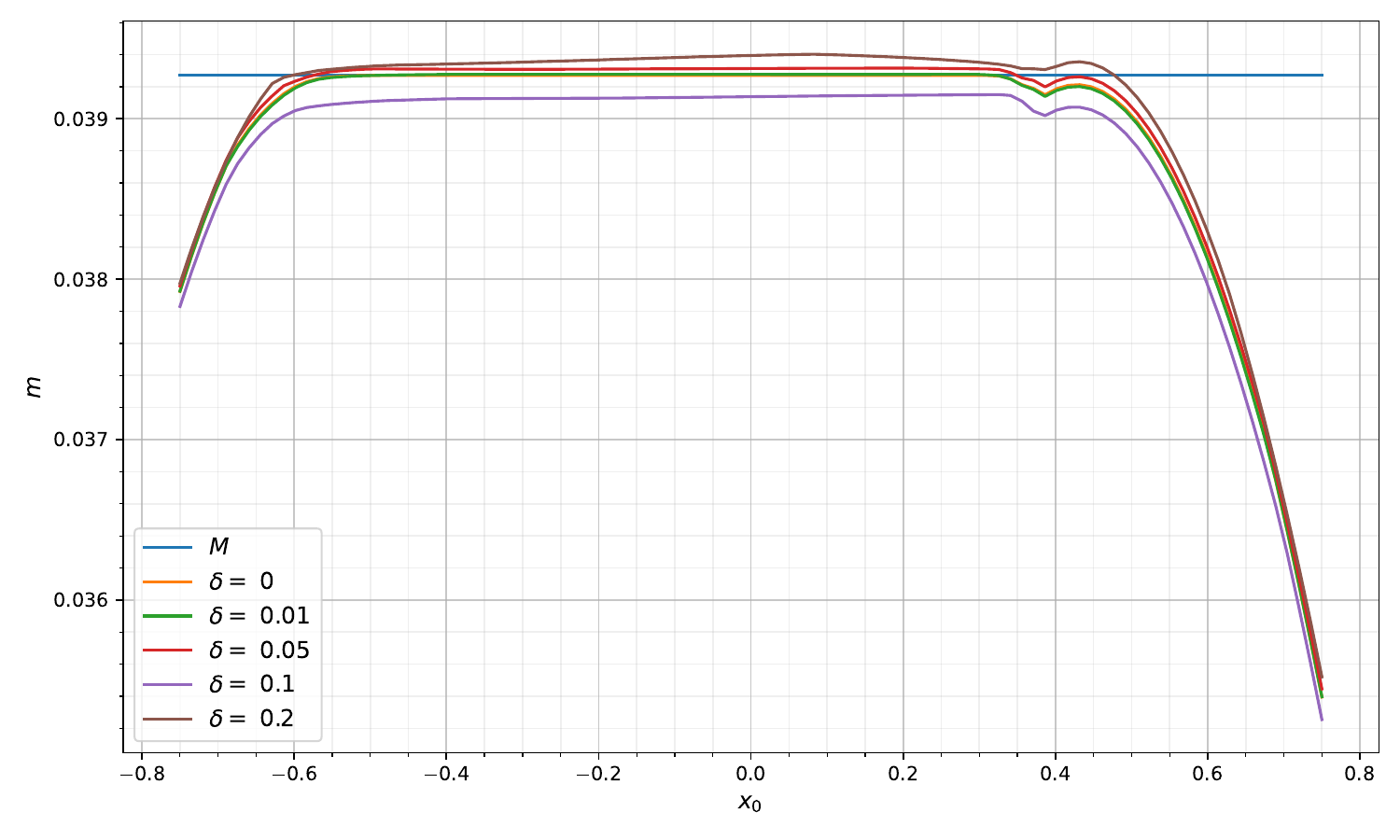}
\caption{Mass for different positions of $S$.}
\label{f-25}
\end{figure}

\begin{figure}[ht]
\centering
\includegraphics[width=0.75\textwidth]{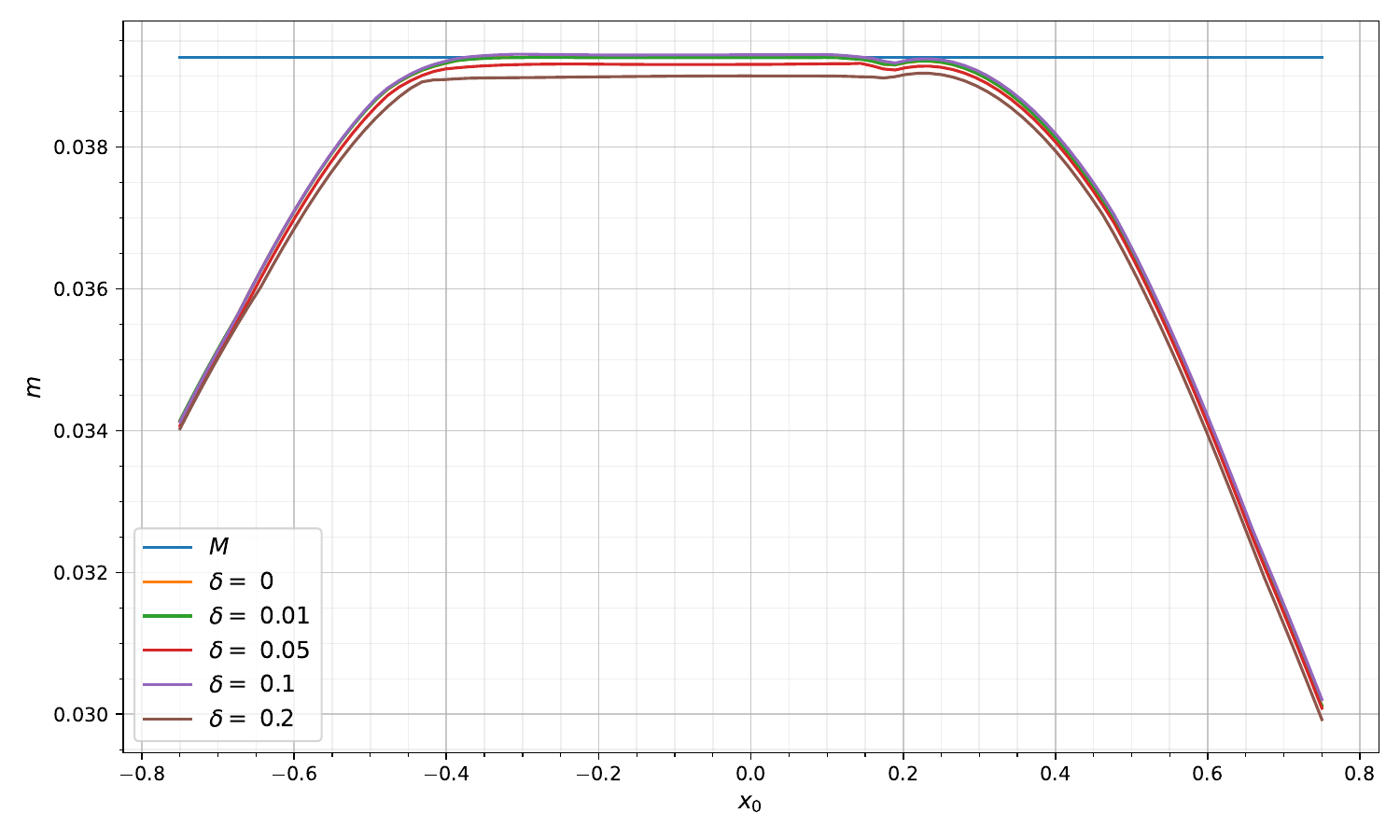}
\caption{Mass for different positions of the reduced window $S$.}
\label{f-26}
\end{figure}

Another practical metric is the total source mass within $S$:
\[
 m = \int_S \mu(\bm y) d\bm y.
\]
For the test problem, the true source mass is $M = 0.0125 \pi$.
Figures~\ref{f-25}, \ref{f-26} show that the maximum mass correlates with both the magnitude and location of sources.
Reducing the window size improves localization accuracy.

\section{Conclusions}

\begin{enumerate}[(1)]
    \item  
The inverse potential problem for the two-dimensional Laplace operator is considered with measurements on a contour that encloses unknown sources. 
In the general case, the solution to this problem is not unique, as the distributed density of the volumetric potential cannot be determined unambiguously.
The inverse problem is formulated to identify the domain containing the sources.
    \item 
The measured volumetric potential is approximated by a single-layer potential on the boundary of the domain that includes the sources.
For an approximate solution of the corresponding first-kind Fredholm integral equation, the classical Tikhonov regularization method is used. It provides a stable approximate normal solution (minimal in norm).
    \item 
In this work, the inverse problem of finding the equivalent single-layer potential is considered under the a priori constraint of nonnegative potential density. 
Such assumptions hold when the volumetric potential density is nonnegative and all sources lie within the support of the single-layer potential. 
After discretization, the approximate solution is determined by minimizing the residual within the class of nonnegative solutions using the Nonnegative Least Squares (NNLS) algorithm.
    \item 
A heuristic algorithm for identifying the support of the volumetric potential (Nonnegative Density Domain, NNDD, algorithm) is proposed, based on approximate computation of single-layer potentials with various supports in the class of positive densities.
Identification is performed either by the residual of the potential on the observation surface or by the total mass of the layer.
    \item 
Numerical experiments are carried out for a test two-dimensional problem with analytically prescribed potential on the observation surface.
For the numerical computation of the single-layer potential, the Tikhonov regularization method and the Nonnegative Density Domain algorithm are applied. 
Calculations with different positions of the single-layer potential support demonstrate the potential capability of identifying the support of the volumetric potential within the class of positive densities.
\end{enumerate}

\begin{funding}
This research was funded by the Committee of Science of the Ministry of Science and Higher Education of the Republic of Kazakhstan (Grant no. BR27100483).
\end{funding}

\providecommand{\bysame}{\leavevmode\hbox to3em{\hrulefill}\thinspace}
\providecommand{\MR}{\relax\ifhmode\unskip\space\fi MR }
\providecommand{\MRhref}[2]{%
  \href{http://www.ams.org/mathscinet-getitem?mr=#1}{#2}
}
\providecommand{\href}[2]{#2}

\end{document}